\documentclass[12pt]{amsart} 

\usepackage{graphicx,amssymb,colordvi,textcomp,latexsym} 

\usepackage{tikz}
\usepackage{ulem}

\usepackage[title]{appendix}
\usepackage[showframe=false]{geometry}
\usepackage{changepage}

\usepackage{tikz}

\usepackage{color}
\usepackage{colortbl}
\usepackage{soul}

\usepackage{array}
\usepackage[hidelinks]{hyperref}
\usepackage{amsfonts}
\usepackage{amsmath}
\usepackage{amsthm}

\usepackage{adjustbox}
\usetikzlibrary{arrows,shapes, positioning, matrix, patterns, decorations.pathmorphing}

\usepackage{color}
\usepackage{colortbl}
\usepackage{soul}

\newcommand{\R}{\mbox{$\mathbb{R}$}}

\newcommand{\ZZ}{\mathbf{Z}}

\newcommand{\NN}{\mathbf{N}}

\newtheorem{lemma}{Lemma}[section]

\newtheorem{prop}[lemma]{Proposition}
\newtheorem{thm}[lemma]{Theorem}

\theoremstyle{definition}
\newtheorem{Def}[lemma]{Definition}
\newtheorem{exam}[lemma]{Example}

\theoremstyle{remark}
\newtheorem{rem}[lemma]{Remark}

\newcommand{\act}{\mathrm{act}}
\newcommand{\rep}{\mathrm{rep}}

\newcommand{\etal}{{\it et al.}}

\usepackage{xcolor}

\title{Synchrony in Gene Regulatory Networks}

\author{Manuela Aguiar}
\address{Manuela Aguiar, Faculdade de Economia, Centro de Matem\'atica, Universidade do Porto,
Rua Dr Roberto Frias, 4200-464 Porto, Portugal.}
\email{maguiar@fep.up.pt}
\author{Ana Dias}
\address{Ana Dias, Departamento de Matem\'atica, Centro de Matem\'atica, Universidade do Porto,
Rua do Campo Alegre, 687, 4169-007 Porto, Portugal}
\email{apdias@fc.up.pt}
\author{Haibo Ruan}
\address{Haibo Ruan, Institute of Mathematics, Technical University of Hamburg, Am Schwarzenberg-Campus 3,
	D-21073 Hamburg, Germany}
\email{haibo.ruan@tuhh.de}

\date{\today} 

\begin{document}

\begin{abstract} 
We consider two mathematical dynamical models of gene regulatory networks (GRNs) and obtain results on robust synchronization on these dynamical models based on the existing theoretical results in the coupled cell network formalism. We also explore the concepts of quotient networks and network lifting in the context of GRNs which are related to the process of gene duplication and the phenomenon of genetic redundancy.
\end{abstract}

\maketitle

\tableofcontents

\section{Introduction and Motivation} \label{sec:intro}

Gene regulatory networks (GRNs) are biochemically interacting systems that are composed of molecular regulators in a cell which interact with each other directly or indirectly through their RNA and protein expression products (\cite{LLC_2019,VCS_2013}). They underlie genetic regulatory mechanisms that determine and control cellular functions such as cell operations, cell-cycle progression and responses to environmental signals (\cite{MSA_2004}).	
 With the continuous progress in genome sequencing technology, an enormous amount of experimental data on gene expression and regulation have been made available, which may provide the first step towards understanding how cells survive, reproduce and adapt (\cite{RRB_2012}). However, given the structural and dynamical interconnectedness of genomes, many cellular processes  that control the development and adaptation of multicellular life forms still remain widely undetermined (\cite{BT_2004,SRG_2008}). 


Mathematical models that offer tractable analysis of dynamical processes in GRNs and their collective behaviour are sought after for understanding, predicting and controlling regulatory mechanisms, see for example \cite{A06, KS08, KLWK16}. The theory of coupled cell networks and their associated dynamical systems, proposed by Golubitsky and Stewart~\cite{GS06}, and by Field~\cite{F04}, have provided mathematical frameworks for understanding collective dynamics on coupled networks such as synchronization and synchrony-related bifurcations, \cite{SGP03, F04, GST05, GS06}. 
Coupled cell systems are dynamical systems ({\it cells}) that are coupled together through mutual interactions, abstracted by the associated network, which exert influences on  the temporal evolution of each other. 
A key advantage of these formalisms is that they allow theoretical deduction of collective dynamics based only on the underlying network structure, without referring to the specifics of each cell. This leads to theoretical results that are relatively independent of modelling specifics of individual cells and are more importantly, naturally compatible with structure-related dynamical processes   in GRNs.

In this paper, we consider two mathematical dynamical models of GRNs and identify their robust patterns of synchrony using theoretical results  from the coupled cell network formalism.  In the context of gene regulations, the synchronization of specific subsets of genes can be a result of gene duplication which also contributes to the phenomenon of genetic redundancy, where two or more genes perform the same function and consequently, the inactivation of one of the genes has little or no effect on the biological phenotype. See, for example, Nowak~\etal~\cite{NBCS1997}.

Moronea, Leifera and Maksea~\cite{MLM2020} introduce the use of network fibration symmetries in the analysis of the transcriptional regulatory network of bacterium {\it Escherichia coli}. They emphasize that {\it fibration symmetries} in GRNs may identify {\it clusters of synchronized genes} which can be collapsed, each group into a single gene.  
In Aguiar, Dias, Golubitsky and Leite~\cite{ADGL09} and its extension to weighted networks by Aguiar and Dias~\cite{AD18}, it is shown that bigger networks can lead to smaller networks called {\it quotients} and where, from the network dynamics point of view, the quotient network dynamics is in correspondence to the bigger network dynamics where genes in each cluster are synchronized.

It is known that if a coupled cell network has permutation symmetry $\Gamma$ (by permuting  the cells) and $\Sigma$ is a subgroup of $\Gamma$, then the fixed point subspace of $\Sigma$ is a flow-invariant subspace for any coupled cell system consistent with the network structure. Moreover, these fixed-point subspaces are described in terms of equalities of cell coordinates. The flow-invariance of a fixed-point subspace under the dynamics of a coupled cell system implies that any solution that starts in that space remains there for all time. That is, a solution where certain cells are synchronized at one time remain synchronized for all time. Furthermore, this phenomena is robust in the sense that, symmetry-preserving perturbations of the systems maintain the flow-invariance of these fixed-point subspaces. Such spaces are called {\it network synchrony spaces}. In particular, this concept applies to GRNs of genes where the underlying networks have symmetry $\Gamma$. That is, for any model equations chosen to describe a GRN, it is expected that the symmetry of the underlying network will be reflected at the equations and the above result will apply. See,  Antoneli and Stewart~\cite{AS06} and Moronea, Leifera and Maksea~\cite{MLM2020}.  As pointed out by Antoneli and Stewart~\cite{AS06}, for general networks there can be synchrony subspaces that are not forced by the symmetries of the graph, if any.
	
 As mentioned earlier, an important aspect of the identification of patterns of synchrony is that, associated with each synchrony pattern, there is a smaller (quotient) network of the GRN, whose dynamics  is related with the total GRN when restricted to the corresponding synchrony subspace. 
This is in accordance with a common technique used in science,  from both theoretical and experimental points of view: to first investigate small networks either from  dynamical or statistical point of view and then to expand the results to larger networks that can be related  to the smaller ones. 
There are different approaches as to how smaller networks can be related to larger ones. 
In this paper, we are following the dynamical perspective and the process of expanding a (smaller) quotient network to a bigger network which ensures that the dynamics associated with the smaller network also occurs at the bigger network. See Aguiar, Dias, Golubitsky and Leite~\cite{ADGL09}. Another perspective that is common to follow is the recognition of network {\it motifs}, that is, small subnetworks that occur in complex networks with frequency  significantly higher than those in randomized networks. See for example Milo, Shen-Orr, Itzkovitz, Kashtan, Chklovskii and Alon~\cite{MSIKCA02}.
 
Dewey and Galas~\cite{DG} discuss a network growth model based on {\it gene duplication} pointing out evidence through a variety of genoma analysis that new genes are almost always created by duplication.  On the network community, we recall work on network dynamics where {\it lifting} is a process of embedding smaller network dynamics into  bigger network dynamics. The terminology used is that of lifting (inflating) cells in a small network so that some of the cells give rise to copies of themselves and where the interactions linking the new cells in the bigger network guarante that  the dynamics of the bigger systems contains the dynamics of the smaller systems. See Aguiar, Dias, Golubitsky and Leite~\cite{ADGL09}, its extension to weighted networks by Aguiar and Dias~\cite{AD18} and Ashwin, Aguiar, Dias and Field~\cite{AADF11}. 
A trivial observation is that, in general, there is no uniqueness in this lifting process. That is, fixing a $k$-cell network, if $n>k$ then there are many $n$-cell networks which are lifts (or inflations) of the $k$-cell network. Moreover, there is a method of enumerating the lifts if the $k$-cell network  is fixed. See Theorem 2.5 of Aguiar, Dias, Golubitsky and Leite~\cite{ADGL09} valid for coupled cell systems following the formalisms of Stewart, Golubitsky and Pivato~\cite{SGP03}, Golubitsky, Stewart and T\"{o}r\"{o}k~\cite{GST05} and Field~\cite{F04}, where the edges in the graph network structure have assigned nonnegative integer numbers. See also the extension of this result to weighted networks, where  the connections have attached strengths that can be any real number and  the coupled cell systems have additive input structure, in  Theorem 2.13 of Aguiar and Dias~\cite{AD18}. For both setups, the enumeration method goes through the characterization of the network adjacency matrices of the larger networks determined by the smaller network adjacency matrices. A nice observation in the enumeration method, is that, fixing $n$, the number of $n$-cell lifts is finite for nonnegative integer matrices and nonfinite in the weighted setup.

In this work we aim to connect and adapt some of the existent results concerning network robust synchronization to GRNs.  We address mainly two issues. One concerns the existence of synchrony spaces in GRNs and their description; the other is related with robustness of such synchrony patterns. That is, does the answer to the first issue depends on the model equations approach? We address the robustness question to SUM and MULT models. From our results we conclude that, in general, the synchronization patterns can be quite distinct from the SUM and MULT model. Moreover, we see that if the activation and repression functions involved at the GRNs models imply structural relation, unexpected synchrony patterns may occur that are arising due to the specificities of the activation and regression functions.

We consider next an GRN  and exemplify the processes of duplication and synchronization of genes.  

\subsection{Example}\label{subsec:intro}
Consider the three transcriptional repressor system used by Elowitz and Leibler~\cite{EL00} to build a repressilator network, in {\it Escherichia coli}. Note that one of the results obtained in \cite{EL00}  is that, depending on the values of several parameters, such as, the dependence of transcription rate on repressor concentration, the translation rate, and the decay rates of the protein and messenger RNA,  at least two types of solutions are possible: the system may converge towards a stable steady state, or the steady state may become unstable, leading to sustained limit-cycle oscillations. See also \cite{BL_2015}.
  In the equations model used in \cite{EL00}, each variable $x_i = (m_i,p_i)$, for $i=1,2,3$, describes the concentrations of two gene products, the mRNAs and proteins, which vary in continuous time, and their time derivatives are expressed as functions of the variables. Graphically, each node in the network 
 represents a gene. Also, an interaction between  two nodes is  represented by an edge. Moreover, in this example there is only a repression type of interaction with specific positive real weight 
 $\alpha$. See Figure~\ref{fig:exEL00}. The {\it repression weight adjacency matrix} is 
 $$
 W^{-} = \left[ w^{-}_{ij} \right] = 
 \left(
 \begin{array}{l|l|l}
 0         & 0 & \alpha \\
 \hline 
 \alpha & 0 & 0  \\
 \hline 
 0 & \alpha & 0
 \end{array}
 \right), 
 $$
 where if  there is an interaction from gene $j$ to gene $i$, the entry $w_{ij}^{-}$ is nonzero and it denotes the positive weight of the repression type interaction. 
 The three-gene network equations used in \cite{EL00} are: 
 \begin{equation}\label{eq:ex2SUM}
 {\footnotesize
\left\{
\begin{array}{l}
\dot{x}_1 = 
A_{\beta} x_1 + \left[ \alpha \rep (p_3) + \alpha_0\right] 
\left(
\begin{array}{l}
1 \\
0
\end{array}
\right)
\\
\dot{x}_2 = 
A_{\beta} x_2 + \left[ \alpha \rep (p_1) + \alpha_0\right] 
\left(
\begin{array}{l}
1 \\
0
\end{array}
\right)
\\
\dot{x}_3 = 
A_{\beta} x_3 + \left[ \alpha \rep (p_2) + \alpha_0 \right] 
\left(
\begin{array}{l}
1 \\
0
\end{array}
\right)
\end{array}
\right. 
\Leftrightarrow 
\left\{
\begin{array}{l}
\dot{x}_1 = 
A_{\beta} x_1 + \left( 
\frac{\alpha}{1 + p_3^2} + \alpha_0 \right) 
\left(
\begin{array}{l}
1 \\
0
\end{array}
\right)
\\
\dot{x}_2 = 
A_{\beta} x_2 + \left( 
\frac{\alpha}{1 + p_1^2} + \alpha_0 \right) 
\left(
\begin{array}{l}
1 \\
0
\end{array}
\right)
\\
\dot{x}_3 = 
A_{\beta} x_3 + \left( 
\frac{\alpha}{1 + p_2^2} + \alpha_0\right) 
\left(
\begin{array}{l}
1 \\
0
\end{array}
\right)
\end{array}
\right. , 
}
\end{equation}
where
$$
A_{\beta} 
= 
\left(
\begin{array}{rr}
-1 & 0 \\
1 & -\beta
\end{array}
\right) \quad \mbox{ and } \quad \rep (p) = 1 - \frac{p^2}{1 + p^2} = \frac{1}{1 + p^2}\, .
$$
Thus all the three genes have the same internal (linear) dynamics and all the weights are equal to $\alpha$. So in fact the system (\ref{eq:ex2SUM}) has cyclic permutation symmetry ${\bf Z}_3$ on the three nodes which justifies nicely the oscillatory behaviour found in \cite{EL00},  
see for example Theorem~XVII 8.2 of Golubitsky, Stewart and Schaeffer~\cite{GSS88} and Chapter 4 of Golubitsky and Stewart~\cite{GS02}. Also, the stabilty of the equilibrium depends on $\alpha$ and $\beta$. Note that, for each $x_i$ gene  equation, only the $m_i$ mRNA  equation depends directly on the interaction from other gene $j$ through the $p_j$ protein concentration.

\begin{figure}
	\begin{center}
\begin{tikzpicture}
 [scale=.15,auto=left, node distance=1.5cm, 
 ]
\node[fill=white,style={circle,draw}] (n1) at (4,0) {\small{1}};
\node[fill=white,style={circle,draw}] (n2) at (24,0) {\small{2}};
\node[fill=white,style={circle,draw}] (n3) at (14,9) {\small{3}};
 \path 
 (n1) [-|] edge[thick] node [near end, above] {{\tiny $\alpha$}} (n2)
 (n2) [-|] edge[thick] node [near end, above] {{\tiny $\alpha$}} (n3)
 (n3) [-|] edge[thick] node [near end, above] {{\tiny $\alpha$}} (n1)
   ;
\end{tikzpicture} 
\caption{A GRN with three genes. It corresponds to the repressilator composed by a cyclic negative-feedback loop of three repressor genes presented in \cite{EL00}.}
\label{fig:exEL00}
\end{center}
\end{figure}
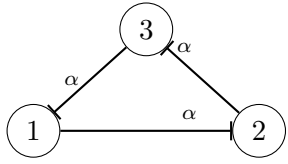

\begin{figure}
	\begin{center}
\begin{tikzpicture}
 [scale=.15,auto=left, node distance=1.5cm, 
 ]
\node[fill=white,style={circle,draw}] (n1) at (4,0) {\small{1}};
\node[fill=white,style={circle,draw}] (n2) at (24,0) {\small{2}};
\node[fill=white,style={circle,draw}] (n3) at (10,9) {\small{3}};
\node[fill=white,style={circle,draw}] (n4) at (18,9) {\small{4}};
 \path 
 (n1) [-|] edge[thick] node [near end, above] {{\tiny $\alpha$}} (n2)
 (n2) [-|] edge[thick] node [near end, above] {{\tiny $\alpha$}} (n3)
 (n2) [-|] edge[thick] node [above] {{\tiny $\alpha$}} (n4)
 (n3) [-|] edge[thick] node [near end, above] {{\tiny $\alpha_1$}} (n1)
 (n4) [-|] edge[thick] node [near end, below] {{\tiny $\alpha_2$}} (n1)
 ;
\end{tikzpicture} 
\caption{A GRN with four genes that can be interpreted from the three gene GRN of Figure~\ref{fig:exEL00} where gene $3$ was duplicated giving rise to genes $3,4$.}
\label{fig:dup}
 \end{center}
\end{figure}
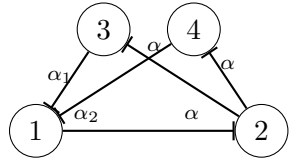

\subsubsection*{Gene duplication example} 
Consider any $4 \times 4$ regression weighted adjacency matrix of the form
$$
 \left(
 \begin{array}{l|l|ll}
 0         & 0 & \alpha_1 & \alpha_2  \\
 \hline 
 \alpha & 0 & 0  & 0 \\
 \hline 
 0 & \alpha & 0 & 0 \\
 0 & \alpha & 0 & 0
 \end{array}
 \right)
 $$
 where the weights $\alpha_1, \alpha_2 > 0$ are such that $\alpha_1 + \alpha_2 = \alpha$, see Figure~\ref{fig:dup}. Take the following model equations for the four gene GRN: 
\begin{equation}\label{eq:exdup}
\left\{
\begin{array}{l}
\dot{x}_1 = 
A_{\beta} x_1 + \left[ \alpha_1 \rep (p_3) + \alpha_2 \rep (p_4) + \alpha_0\right] 
\left(
\begin{array}{l}
1 \\
0
\end{array}
\right)
\\
\dot{x}_2 = 
A_{\beta} x_2 + \left[ \alpha \rep (p_1) + \alpha_0\right] 
\left(
\begin{array}{l}
1 \\
0
\end{array}
\right)
\\
\dot{x}_3 = 
A_{\beta} x_3 + \left[ \alpha \rep (p_2) + \alpha_0 \right] 
\left(
\begin{array}{l}
1 \\
0
\end{array}
\right)\\
\dot{x}_4 = 
A_{\beta} x_4 + \left[ \alpha \rep (p_2) + \alpha_0 \right] 
\left(
\begin{array}{l}
1 \\
0
\end{array}
\right)
\end{array}
\right.  
\end{equation}
where, as above $\rep (p) = 1/(1 + p^2)$.  
Restricting these equations to the space where the variables corresponding to genes $3,4$ are identified, that is, $\{ x:\, x_3 = x_4\}$, we obtain equations (\ref{eq:ex2SUM}). In particular, an oscillatory solution of the system (\ref{eq:ex2SUM}) is in correspondence to an oscillatory solution of system (\ref{eq:exdup}) where genes $3,4$ quantities are synchronized for all time. 

The four-gene GRN of Figure~\ref{fig:dup} is said to be  a {\it lift} network of the three-gene GRN of Figure~\ref{fig:exEL00}. Also, it  is used the terminology that  the three-gene GRN of Figure~\ref{fig:exEL00}  is a {\it quotient} of  the four-gene GRN of Figure~\ref{fig:dup} by the synchrony space $\{ x:\, x_3 = x_4\}$. 

We can interpret the four-gene GRN in Figure~\ref{fig:dup} obtained from the three-gene GRN of Figure~\ref{fig:exEL00} by {\it duplication} of gene 3 into genes $3,4$. 
We can also look at the  three-gene GRN of Figure~\ref{fig:exEL00} as obtained from the four-gene  GRN in Figure~\ref{fig:dup} by {\it collapsing} the $3,4$ gene  cluster into a single gene $3$. This cluster exists due to the fibration symmetries of the GRN of Figure~\ref{fig:dup}.  Equivalently, this is due to the facts that both genes $3,4$ receive a repression type  interaction with the same weight from (the cluster formed by the single) gene $2$, and the sum of weights of the repression type  interactions from the cells in the cluster of $3,4$ to gene $1$ is equal to the weight of the repression type interaction to gene $1$.

The paper is organized in the following way. In Section~\ref{sec:prel} we provide some basics of gene regulatory networks that we follow in this paper. We present two dynamical models,  the SUM and MULT models, which differ  in  the way the GRN regulatory function operates over the gene inputs. In Section~\ref{sec:synonGRN}, we introduce the concept of synchrony pattern (synchrony subspace) for the two types of GRN models equations and characterize the synchrony patterns for GRNs. See Propositions~\ref{prop:synSUM}-\ref{prop:SUMquotyes} and Propositions~\ref{prop:MULTsyn}-\ref{prop:MULT_quo}. {When the regulatory functions are Hill-like (cf. (\ref{eq;reg_Hill_like_0})-(\ref{eq;reg_Hill_like_infty})),} a corollary of these results is that the GRN synchrony patterns are completely determined by the GRN  structural matrices. 
In Section~\ref{sec:enumlift}, we combine Theorem 2.13 of Aguiar and Dias~\cite{AD18} with Propositions~\ref{prop:synSUM}-\ref{prop:SUMquotyes} and Propositions~\ref{prop:MULTsyn}-\ref{prop:MULT_quo}, obtaining a characterization method of the $n$-gene GRNs that are lifts of a fixed $k$-gene GRN, for  the SUM using the  $k$-gene GRN activation and regression weighed matrices and for MULT models, using the  $k$-gene GRN activation and regression weighed matrices, and the multiplicities matrices. See Theorems~\ref{thm:enumSUM}-\ref{thm:enumMULT}.
In Section~\ref{sec:regdepsyn}, we remark that for models where the activation and repression regulatory functions are structurally related, then there may be other synchrony patterns forced by those regulator functions, which we call regulatory dependent synchrony patterns.  
We give an example of such regulatory-dependent synchrony patterns in the context of circadian rhythmic clock models in Section \ref{sec:mammalian}. Finally, in Section~\ref{sec:future}, 
we summarize the work presented and point out directions for future work.

\section{Preliminaries}\label{sec:prel}

\subsection{Gene Regulatory Networks (GRNs)}

Gene regulatory networks can be modelled by two classes of mathematical models: the Boolean model (or discrete model) and the differential equation model (or continuous model). In boolean models, the activity of each gene is expressed in one of two states, ON or OFF, and the state of a gene is determined by a Boolean function of the states of related genes. In the differential equation models, the variables describe the concentrations of gene products such as mRNAs and proteins as continuous values, and their time derivative is expressed as a function of the variables themselves (cf. \cite{Chesi2010}).

In this paper, we follow the differential equation model.  The activity of a gene is regulated by other genes through the concentrations of their gene products, which function as transcription factors. Regulation can be quantified by the ``response characteristics'' given by the level of gene expression as  a function of transcription factors (cf. \cite{LCA_2006}). Formally, 
a GRN  can be described by a set of ordinary differential equations, for $i=1,2,\dots,n$,  of the form
\begin{equation}\label{eq:GRN}
\begin{cases}
\dot m_i(t)=-a_i m_i(t)+b_i(p_1(t),p_2(t),\dots, p_n(t))\\
\dot p_i(t)=-c_ip_i(t)+d_im_i(t)
\end{cases}, 
\end{equation}
where $m_i,p_i\in \R$ are the concentrations of mRNA and protein, respectively, of the $i$-th node, $a_i,c_i>0$ are degradation rates of mRNA and 
protein, respectively,  and $d_i>0$ is a constant.

Graphically, the set of equations (\ref{eq:GRN}),  for $i=1,2,\dots,n$, translates a GRN as a network  of $n$ nodes, where each node $x_i:=(m_i,p_i)$ is connected with {another node if there is} an interaction {between them} via the regulatory function $b_i$.

\subsection{ Models of GRNs and Regulatory Functions}
 The regulatory function $b_i$ in (\ref{eq:GRN}) plays a key role in the dynamical modelling of GRNs. It is generally a nonlinear function of its variables $p_1(t), p_2(t),\dots,  p_n(t)$ and may assume monoticity in each variable in simplified cases  (cf. \cite{T_2010}, more references from \cite{LCA_2006}). Depending on all biochemical reactions involved, it may have a complicated form which is,  in practice, determined heuristically.
	
There are typically two regulatory logics behind $b_i$, depending on whether each transcription factor acts {\it additively} or {\it multiplicatively} to regulate the $i$-th gene. It is called {\it SUM logic}, if
\[b_i=\sum_{j} b_{ij}(p_j(t))\]	 
and it is called {{\it  MULT} logic}, if  
\[b_i=\Pi_{j} b_{ij}(p_j(t)).\]
See for example Chesi~\cite{Chesi2010}, Chesi and Hung~\cite{CH08} and \cite{KKBLRH_2014,PKBKH_2018,PKMKH_2016}.

There are two types of regulatory functions: those that {\it activate} and those that {\it repress}  the target gene expression in the $i$-th node,  described by a monotonically {\it increasing} and {\it decreasing} function $b_i$, respectively.  

 
One of the most commonly used {type of activation} regulatory function is described by (cf. ~\cite{Elowitz_1999,Chesi2010} for example) the {\it Hill function}
 \begin{equation}\label{eq:Hill_f}
 \act  (p) = \frac{p^n}{\beta^n + p^n}, 
  \end{equation}
 where  $n\in \NN$ characterizes the steepness of regulation and $\beta$  marks the mid-value of maximal reachable value of $\act$. See Remark \ref{rem:GRF} for details. The repression is {frequently} modelled by
 \begin{equation}\label{eq:Hill}
 \rep (p)=1-\act(p) = \frac{\beta^n}{\beta^n + p^n}\, .
 \end{equation}

For our purpose, to allow more general modelling possibilities, we will only assume that ${\mathrm act}$ and ${\mathrm rep}$ are strictly monotonic between 0 and 1, and satisfy
\begin{align}
\act(p)\to 0  &\ \mbox{ and } \ \rep (p)\to 1 \  \mbox{ when } \  p\to 0 \label{eq;reg_Hill_like_0}\\
\act (p)\to 1  &\ \mbox{ and }  \  \rep (p)\to 0\ \mbox{ when } \  p\to \infty, \label{eq;reg_Hill_like_infty} 
\end{align}
 and call them the {\it Hill-like} regulatory functions.

Consider a GRN of $n$ nodes (genes), $x_i= (m_i,p_i) \in \R^{+}_{0}$ for $i= 1, \ldots, n$, where the concentration $m_i$ of mRNA and $p_i$ of protein are measured for each node. Take the internal dynamics function of the $i$th  node {as} the $2 \times 2$ matrix 
$$
A_i = 
\left( 
\begin{array}{rr}
-a_i & 0 \\
d_i & -c_i
\end{array}
\right),
$$
where $a_i, c_i, d_i$ are positive real constants.  That is, we assume that the internal dynamics of the genes is linear.  In particular, it follows that each $A_i$ is invertible. 
Given $i$, define 
$$
\begin{array}{l} 
I_i^+ = \{ j \in \{ 1, \ldots, n\}:\, \mbox{ gene } j \mbox{ activates gene } i\}, \\
I_i^- = \{ j \in \{ 1, \ldots, n\}:\, \mbox{ gene } j \mbox{ represses gene } i\}
\end{array}.
$$

\begin{rem}
We allow $I_i^+\cap I_i^-\ne \emptyset$. That is, a gene can activate and repress a same gene.
\hfill $\Diamond$
\end{rem}

Denote the two $n \times n$  network adjacency matrices by $W^+ = [w^+_{ij}]$ and  $W^- = [w^{-}_{ij}]$. If there is activation (resp. repression) from gene $j$ to gene $i$, we have that $w^+_{ij} > 0$ (resp. $w^{-}_{ij} > 0$) stands  for the maximal achievable level of the activation (resp. repression)  from gene $j$ to gene $i$. Otherwise, $w^+_{ij} =0$ (resp. $w^{-}_{ij} = 0$). We call $W^+$ (resp. $W^{-}$) the network {\it weighted activation} (resp. {\it repression) adjacency matrix}.

\subsection*{The SUM model}

In the SUM model, the regulatory function is operating by {\it addition}:
\begin{equation}\label{eq:SUM}
\dot{x}_i = A_i x_i + 
\sum_{j \in I^{-}_i} w^{-}_{ij} \rep (p_j) 
\left(
\begin{array}{l}
1 \\
0
\end{array}
\right)
+ 
\sum_{j \in I^{+}_i} w^+_{ij} \act  (p_j) 
\left(
\begin{array}{l}
1 \\
0
\end{array}
\right) \qquad \left( i=1, \ldots, n\right)\, .
\end{equation}

\subsection*{The MULT model}

In the MULT model, the regulatory function is operating by {\it multiplication}:
\begin{equation}\label{eq:MULT0}
\dot{x}_i = A_i x_i + 
\prod_{j \in I^{-}_i} w^{-}_{ij} \rep^{m_{ij}^-} (p_j) \, 
\prod_{j \in I^{+}_i} w^{+}_{ij} \act^{m_{ij}^+} (p_j) 
\left(
\begin{array}{l}
1 \\
0
\end{array}
\right)  \qquad \left( i=1, \ldots, n\right), \, 
\end{equation}
where $m_{i,j}^{\pm}$ indicates the multiplicative regulation of the gene $i$ from the gene $j$. Graphically, these multiplicities are usually represented by multiple arrows. However, in the case of $0-1$ adjacency matrices we frequently indicate the multiplicity of the arrows with a number next to each arrow. 
In the general case, if we have both weights $w_{ij}^{\pm}$ and exponents $m_{ij}^{\pm}$ in the model, then we reserve the number next to the arrow for $w_{ij}^{\pm}$ and use multiple arrows to indicate $m_{ij}^{\pm}$. 
Let $M^\pm=(m_{ij}^\pm)$ 
be the {\it activation and repression multiplicity matrices} with the assumption that each entry of $W^{+}$ (resp. $W^{-}$) is nonzero if and only if each entry of $M ^{+}$ (resp. $M^{-}$) is nonzero. The weights $w_{ij}^\pm$, on the other hand, can be combined into one weight 
$$w_i= \prod_{j \in I^{-}_i} w^{-}_{ij}\prod_{j \in I^+_i} w^+_{ij}  \qquad \left( i=1, \ldots, n\right) \, .$$
 Thus, (\ref{eq:MULT0}) is effectively
\begin{equation}\label{eq:MULT}
\dot{x}_i = A_i x_i + w_i
\prod_{j \in I^{-}_i} \rep^{m_{ij}^-} (p_j) \, 
\prod_{j \in I^{+}_i} \act^{m_{ij}^+} (p_j) 
\left(
\begin{array}{l}
1 \\
0
\end{array}
\right)  \qquad \left( i=1, \ldots, n\right) \, . 
\end{equation}


\begin{rem}
In some literature, see for example~\cite{Chesi2010}, it is used a PROD model, which corresponds to a particular case of the MULT model where $m_{ij}^{\pm} = 1$, for all $i,j$. That is, the PROD model equations are given by:
\begin{equation}\label{eq:PROD}
\dot{x}_i = A_i x_i + 
\prod_{j \in I^{-}_i} w^{-}_{ij} \rep (p_j) \, 
\prod_{j \in I^{+}_i} w^{+}_{ij} \act (p_j) 
\left(
\begin{array}{l}
1 \\
0
\end{array}
\right)  \qquad \left( i=1, \ldots, n\right) \, .
\end{equation}
As we will see later, a difficulty we found at the PROD model is that the restriction of PROD model equations to a synchrony 
space is not necessarily a PROD model equations.
\hfill $\Diamond$
\end{rem}

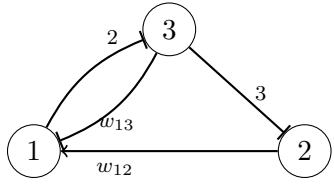
\begin{figure}
	\begin{center}
\begin{tikzpicture}
 [scale=.18,auto=left, node distance=1.5cm, 
 ]
\node[fill=white,style={circle,draw}] (n1) at (4,0) {\small{1}};
\node[fill=white,style={circle,draw}] (n2) at (24,0) {\small{2}};
\node[fill=white,style={circle,draw}] (n3) at (14,9) {\small{3}};
 \path 
 (n1) [-|] edge[bend left=20,thick] node [near end, above] {{\tiny $2$}} (n3)
 (n2) [->] edge[thick] node [near end, below] {{\tiny $w_{12}$}} (n1)
 (n3) [-|] edge[bend left=20,thick] node [below] {{\tiny $w_{13}$}} (n1)
  (n3) [-|] edge[thick] node [near end, above] {{\tiny $3$}} (n2)
   ;
\end{tikzpicture}  
\caption{A GRN with three genes.}
\label{fig:ex2Chesi}
	\end{center}
\end{figure}

\begin{exam} \label{ex:Chesi_ex_2} 
The MULT (PROD) model equations in Example 2 of  \cite{Chesi2010}  correspond to the 
three-cell GRN which appears in Figure\ \ref{fig:ex2Chesi} and have repression and activation adjacency and multiplicities matrices given by 
$$
{\tiny 
W^{-} = \left( 
\begin{array}{lll}
0 & 0 & w_{13}\\
0 & 0 & 3 \\
2 & 0 & 0
\end{array}
\right), \ 
W^{+} = \left( 
\begin{array}{lll}
0 & w_{12} & 0\\
0 & 0 &0 \\
0 & 0 & 0
\end{array}
\right), \ 
M^{-} = \left( 
\begin{array}{lll}
0 & 0 & 1\\
0 & 0 & 1 \\
1 & 0 & 0
\end{array}
\right), \ 
M^{+} = \left( 
\begin{array}{lll}
0 & 1& 0\\
0 & 0 &0 \\
0 & 0 & 0
\end{array}
\right), }
$$
where $w_{12} w_{13} = 4$, respectively. Note that gene $1$ represses gene $3$, gene $2$ activates gene $1$, and gene $3$ represses genes $1, 2$. Moreover, 
all activations and repressions have multiplicity one. The activation function is the Hill function
$\act (p) = p^2/(1 + p^2)$ and the repression function is $\rep (p) = 1/(1 + p^2)$. The  MULT (PROD) model equations 
in Example 2 of  \cite{Chesi2010}  are given by  
\begin{equation}\label{eq:ex2PROD}
\left\{
\begin{array}{l}
\dot{x}_1 = 
A_1 x_1 + w_{12} \act (p_2)  w_{13} \rep (p_3) 
\left(
\begin{array}{l}
1 \\
0
\end{array}
\right)
\\
\dot{x}_2 = 
A_2 x_2 + 3 \rep (p_3) 
\left(
\begin{array}{l}
1 \\
0
\end{array}
\right)
\\
\dot{x}_3 = 
A_3 x_3 + 
2 \rep (p_1) 
\left(
\begin{array}{l}
1 \\
0
\end{array}
\right)
\end{array}
\right.
\Leftrightarrow 
\left\{
\begin{array}{l}
\dot{x}_1 = 
A_1 x_1 + 
\frac{4 p_2^2}{(1+p_2^2) (1 + p_3^2)} 
\left(
\begin{array}{l}
1 \\
0
\end{array}
\right)
\\
\dot{x}_2 = 
A_2 x_2 + 
\frac{3}{1 + p_3^2} 
\left(
\begin{array}{l}
1 \\
0
\end{array}
\right)
\\
\dot{x}_3 = 
A_3 x_3 + 
\frac{2}{1 + p_1^2} 
\left(
\begin{array}{l}
1 \\
0
\end{array}
\right)
\end{array}
\right. ,
\end{equation}
where the internal gene dynamics is linear and determined by certain $2 \times 2$ matrices $A_1, A_2, A_3$.
\hfill $\Diamond$
\end{exam}

\subsection{An example of mammalian circadian rhythm}\label{sec:mammalian}

Following a statistical mechanical framework proposed in \cite{Numbers_2005, Numbers_2005b}, 
a concise GRN model was proposed in \cite{KKBLRH_2014,PKMKH_2016}, which consists of only $5$ gene variables  {\it Bmal1,  Rev-erb-$\alpha$, Per2, Cry1, Dbp} in studying the complex gene regulatory dynamics of the mammalian circadian oscillator. It was shown to describe the known phase relations, amplitudes and wave forms of clock gene expression profiles. Based on the same $5$-gene model, different feedback loops (as sub-networks of the $5$-gene network) have been identified to fit circadian gene expression profiles for different mammalian tissues (\cite{PKBKH_2018}). 

The regulator functions are modelled by (cf. \cite{KKBLRH_2014}, Supplement)
\begin{align}
& \rep(x,b)=\frac{1}{1+bx}\label{eq:5_rep}\\
& \act(x,a)=\frac{1+ax}{1+x}\label{eq:5_act},
\end{align} 
where $a\ge 1$ is a parameter for fold activation and $b\ge 0$ is a repression parameter.  We compare them with the classic Hill regulatory functions in regulation modelling parameters.

\begin{rem}\label{rem:GRF} Many gene regulation functions have been proposed heuristically. Besides the monotonicity with finite asymptotic values, the key features are usually characterized by steepness of regulation, basal transcription or transcription initiation (detected in the absence of transcription factors), the required transcription factor level and the maximal achievable level. 
		
We explain using activators as an example. In case of the regulation given by the Hill functions (cf.  (\ref{eq:Hill_f})-(\ref{eq:Hill})), the transcription fold change is described by
$$H(x)=\frac{a x^n}{b^n+x^n},$$
where $a>1$ is the maximal achievable transcription factor level, $b$ is where the mid-value $\frac a2$ is attained and $n\in \NN$ signifies  the steepness of the curve measured by the slope at the mid-point
\begin{equation}\label{eq:H_steepness}
H'(b)=\frac{a}{4b}n.
\end{equation}

In case of the activating function used for circadian clock models, the transcription fold change is described by
$$C(x)= \displaystyle \left(\frac{1+a\cdot\frac xb}{1+\frac xb}\right)^n=\left(\frac{b+ax}{b+x}\right)^n,$$
where $\frac xb$ is the normalized concentration,  $a^n>1$ is the maximal  achievable transcription factor level, $b$ is where the ``mid-value'' $\big(\frac {a+1}{2}\big)^n$ is attained and $n\in \NN$ signifies  the steepness of the curve measured by the slope at the mid-point
\begin{equation}\label{eq:C_steepness}
C'(b)= \displaystyle \frac{a-1}{4b}\left(\frac {a+1}{2}\right)^{n-1}n\, .
\end{equation}
See Figure \ref{F:GRF_compare}.
 
	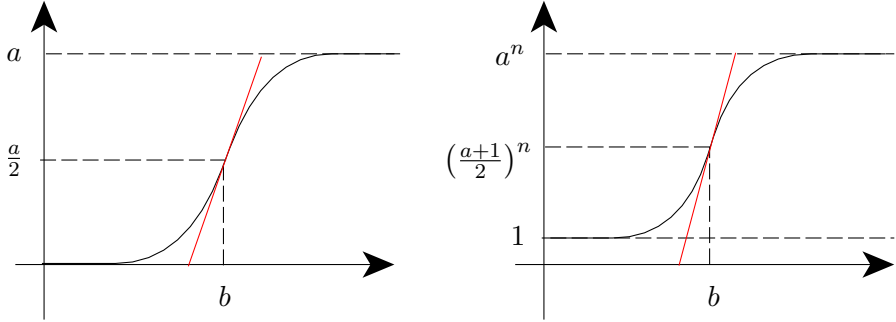
\begin{figure}[htb!]
\begin{tikzpicture}
	\node (n2) at (0,0) {	\includegraphics[width=.9\textwidth]{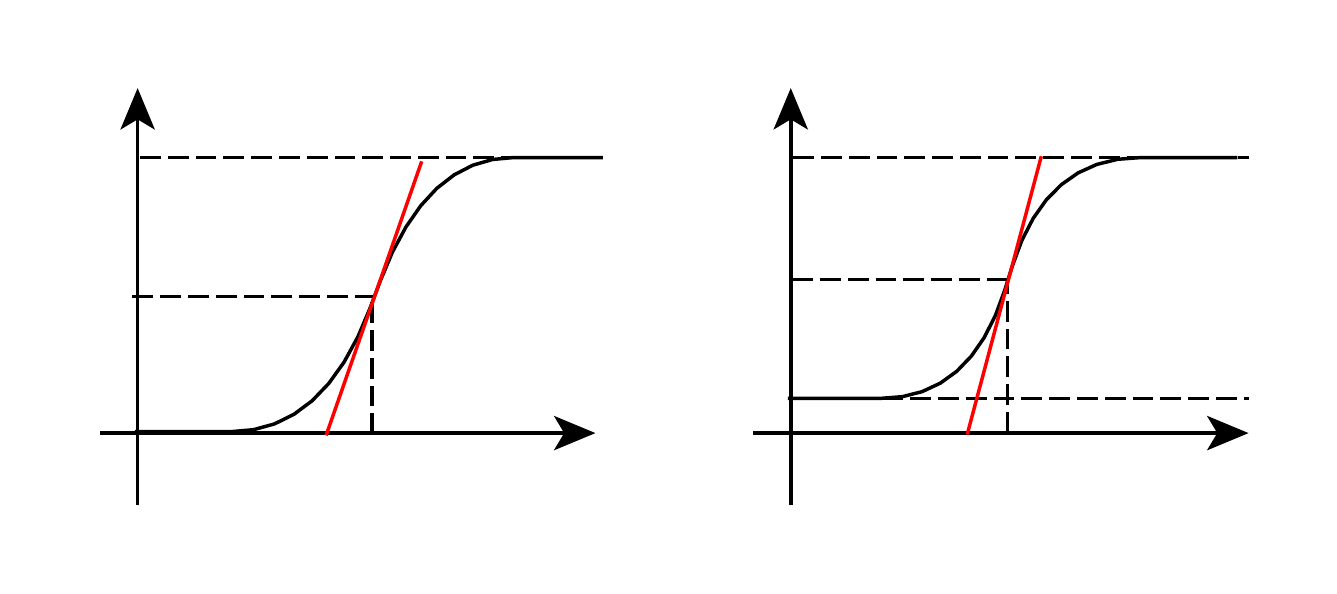}};
\node[fill=white] at (-5.8,1.4) {\small{$a$}};
\node[fill=white] at (.8,1.4) {\small{$a^n$}};
\node[fill=white] at (-5.8,0) {\small{$\frac a2$}};
\node[fill=white] at (.5,0) {\small{$\big(\frac {a+1}{2}\big)^n$}};
\node[fill=white] at (-3,-1.8) {\small{$b$}};
\node[fill=white] at (3.5,-1.8) {\small{$b$}};
\node[fill=white] at (.9,-1) {\small{$1$}};
\end{tikzpicture} 
\caption{Two gene regulation functions for activating transcription factors. The Hill function $H(x)=\frac{a x^n}{b^n+x^n}$ (left) and the activating function $C(x)=\big(\frac{b+ax}{b+x}\big)^n$ used in circadian rhythm model (right). The steepness is measured at the mid-point (cf. (\ref{eq:H_steepness})-(\ref{eq:C_steepness})).}
	\label{F:GRF_compare}  
	\end{figure}

These two functions give very similar qualitative descriptions of gene regulation activities including monotonicity, maximal achievable levels given by $a>1$ and  the steepness factor $n$. 
	
	One crucial difference, however, is the basal transcription given by the asymptotic value at $x=0$, which corresponds to how much transcription can be detected in the absence of any activators. The Hill function $H$ assumes no transcription activities if no activators are present, the function $C$ models with a non-zero constant basal transcription, which in case of circadian models, is modulated by circadian transcription factors and input functions (cf. \cite{KKBLRH_2014, PKBKH_2018, PKBKH_2018}). 
\hfill $\Diamond$
\end{rem}


 \begin{exam} \normalfont 
 The circadian core  $5$-gene clock model is composed of genes  {\it Bmal1,  Rev-erb-$\alpha$, Per2, Cry1,Dbp}, which interact through transcriptional feedback loops of negative and positive regulations. It describes  gene regulatory dynamics of the mammalian circadian oscillator  (cf. \cite{KKBLRH_2014,PKMKH_2016}). See Figure \ref{F:5_gene}. 
\begin{figure}[!htb]
\begin{center}
	\includegraphics[width=.35\textwidth]{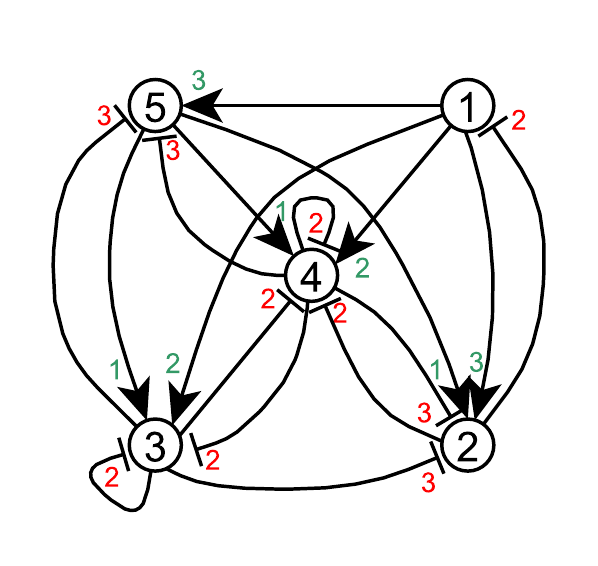} 
	\caption{The circadian core $5$-gene clock model, where cells $1,2,3,4,5$ correspond to the concentration of genes {\it Bmal1,  Rev-erb-$\alpha$, Per2, Cry1,Dbp}, respectively. The number next to the input arrow indicates the number of  binding sites in the target gene (the red for negative, the green for positive feedback).}
	\label{F:5_gene}  
\end{center}
\end{figure}  
The equations employed to describe the dynamics are (cf. \cite{KKBLRH_2014})
\begin{equation} \label{eq:5_gene}
\left\{ 
\begin{array}{ll}
\dot x_1&= \rep^{\scalebox{.7}{\boxed{2}}}(x_2,a_{12})-d_1 x_1\notag\\
\dot x_2&= \act^{\scalebox{.7}{\boxed{3}}}(x_1,a_{21})\act(x_5,a_{25})\rep^{\scalebox{.7}{\boxed{3}}}(x_3,1)\rep^{\scalebox{.7}{\boxed{3}}}(x_4,1)-d_2x_2\notag\\
\dot x_3&= \act^{\scalebox{.7}{\boxed{2}}}(x_1,a_{31})\act(x_5,a_{35})\rep^{\scalebox{.7}{\boxed{2}}}(x_3,1)\rep^{\scalebox{.7}{\boxed{2}}}(x_4,1)-d_3x_3 \\
\dot x_4&= \act^{\scalebox{.7}{\boxed{2}}}(x_1,a_{41})\act(x_5,a_{45})\rep^{\scalebox{.7}{\boxed{2}}}(x_2,1)\rep^{\scalebox{.7}{\boxed{2}}}(x_3,1)\rep^2(x_4,1)-d_4x_4\notag\\
\dot x_5&= \act^{\scalebox{.7}{\boxed{3}}}(x_1,a_{51})\rep^{\scalebox{.7}{\boxed{3}}}(x_3,1)\rep^{\scalebox{.7}{\boxed{3}}}(x_4,1)-d_5x_5,\notag
\end{array} 
\right.
\end{equation}
where $\rep$, $\act$ are given by (\ref{eq:5_rep})-(\ref{eq:5_act}) and the exponents are given by the number of binding sites available in the target gene. The boxed numbers for indicate that these genes $x_j$ for $j\in I_i$ share the same number of binding sites available in the target gene $x_i$. See the numbers in  Figure \ref{F:5_gene}. 
\hfill $\Diamond$
\end{exam}

\begin{rem}\label{rem:mot_osc}
As mentioned above, in  \cite{KKBLRH_2014}, the authors used (\ref{eq:5_gene}) to model circadian clocks, where the boxed numbers are kept equal for genes sharing the same number of binding sites (that is available sites are always occupied) in the target gene. In \cite{PKMKH_2016} as an extension of work from \cite{KKBLRH_2014}, the authors used the same model to run through all oscillating networks and found statistical significance of a sub-network of $\ZZ_3$-symmetry composed of  {\it Rev-erb-$\alpha$, Per2, Cry1} (corresponding to $x_2,x_3,x_4$) that works as a repressilator.
\hfill $\Diamond$ 
\end{rem}

\section{Synchrony in gene regulatory networks}\label{sec:synonGRN}

Motivated by the  formalisms of Stewart, Golubitsky and Pivato~\cite{SGP03}, Golubitsky, Stewart and T\"{o}r\"{o}k~\cite{GST05} and Field~\cite{F04} on coupled cell networks and robust network synchronization for the associated coupled cell systems, we define synchronization partitions for GRN model equations.


\begin{Def} \normalfont 
The {\it gene equivalence partition} of a GRN is the partition  of the genes of the GRN where each part is formed by all the genes with the same internal dynamics function. 
Thus if  $A_{i_1} = A_{i_2}$ (in (\ref{eq:SUM}) or (\ref{eq:MULT})), we have that $i_1,i_2$ belong to the same part of the GRN gene equivalence partition. 
\hfill $\Diamond$
\end{Def}

Consider a  GRN and  an associated dynamical equations model as, for example, the SUM or MULT equations models, (\ref{eq:SUM}) and  (\ref{eq:MULT}), respectively. 

\begin{Def} \normalfont 
Given a partition $P$ of the gene set of the GRN, define the space $\Delta_P$  to be the polydiagonal where gene products concentrations corresponding to the same class of genes  in $P$  are identified. We say that $P$ is a {\it synchronization partition} of the GRN gene set, for the particular model considered, when $\Delta_P$ is flow-invariant under the equations of that model, for any given  Hill-like regulator functions  (cf. (\ref{eq;reg_Hill_like_0})-(\ref{eq;reg_Hill_like_infty})). 
In that case, we call $\Delta_P$ a {\it synchrony pattern} or  a {\it synchrony space} of the GRN. 
\hfill $\Diamond$
\end{Def}

In what follows, given a gene partition $P$ of a GRN, we denote by $[i]$ the part of $P$ that contains gene $i$. Moreover, 
we assume that a necessary condition for genes to synchronize in a robust way 
is that the genes have the same internal dynamics. That is, a synchronization partition $P$ has to refine the gene equivalence partition. 

\begin{Def} \normalfont 
Consider a GRN and an associated dynamical equations model.
Given a synchronization partition $P$ of the gene set of the GRN for the fixed model, if the restriction of the model equations  to the synchrony subspace $\Delta_P$ are equations of the same model type then they are associated to a {\it quotient} GRN. The genes of the quotient GRN correspond to the  parts in the synchronization partition $P$; the activation (repression) interactions between two genes in the quotient are the projection of the activation (repression) interactions between the genes in the corresponding parts of $P$ in the original GRN. We also say that the original GRN is a {\it lift} of the quotient GRN. 
\hfill $\Diamond$
\end{Def}


We address now the issue of synchronization for gene regulatory networks, assuming the models equations are the SUM or the MULT models. More precisely, we characterize the synchrony patterns for GRNs, for the SUM and PROD models, using the activation and repression adjacency and multiplicity matrices. From that, it follows that for the same GRN, distinct patterns of synchrony can occur for the two models, which in particular, lead to distinct dynamical properties for the corresponding dynamical systems. As already remarked, the synchrony subspaces that are forced by the symmetries of a symmetric GRN graph, the fixed-point subspaces,  occur for both SUM and MULT GRN models considered here, as the associated GRNs equations for both models inherit the symmetries of the GRN.  But, there can be synchrony subspaces that are not forced by the symmetries of the graph, if any. See, Antoneli and Stewart~\cite{AS06}. 


\subsection{Synchrony for the SUM model}

\begin{prop}\label{prop:synSUM}
Take an $n$-gene GRN with adjacency matrices $W^+$ and $W^{-}$ and a partition $P$ of the gene set into classes $C_1, \ldots, C_m$ refining the gene equivalence class. The partition $P$ corresponds to a synchrony pattern $\Delta_P$ for the SUM equations model (\ref{eq:SUM}) if and only if for each part $C$, we have that for $k=1, \ldots, m$,  
\begin{equation}\label{cond:synSUM+}
		\sum_{j \in I_i^+ \cap C_k} w^+_{ij} \mbox{ is constant for } i \in C
		\end{equation}
and
\begin{equation}\label{cond:synSUM-}
\sum_{j \in I_i^{-} \cap C_k} w^{-}_{ij} \mbox{ is constant for } i \in C\, .
		\end{equation}
\end{prop}

\begin{proof} Assume that (\ref{cond:synSUM+})-(\ref{cond:synSUM-}) hold for $k=1,\dots,m$ on any equivalence class $C$ under partition $P$. Then, we have for any $i\in C$ that
	\begin{align*}
	\dot x_i&=A x_i+\sum_{j \in I^{-}_{i}} w^{-}_{ij} \rep (p_j) 
	\left(
	\begin{array}{l}
	1 \\
	0
	\end{array}
	\right)
	+ 
	\sum_{j \in I^{+}_
		{i}} w^+_{ij} \act (p_j) 
	\left(
	\begin{array}{l}
	1 \\
	0
	\end{array}
	\right)\\
	&=A x_i+\sum_{k=1}^m \underbrace{\sum_{j \in I^{-}_{i}\cap C_k} w^{-}_{ij}}_{w^-_{i,k}} \rep (p_j) 
	\left(
	\begin{array}{l}
	1 \\
	0
	\end{array}
	\right)
	+ 
	\sum_{k=1}^m \underbrace{\sum_{j \in I^{+}_
		{i}\cap C_k} w^+_{ij}}_{w^+_{i,k}} \act (p_j) 
	\left(
	\begin{array}{l}
	1 \\
	0
	\end{array}
	\right),
	\end{align*}
	where the sums $w^\pm_{i,k}$ by (\ref{cond:synSUM+})-(\ref{cond:synSUM-}), remain the same in the equivalence class of $i$. Thus, for any $i\in C$, the dynamics of $x_{[i]}\in \Delta_P$ is governed by 
	\begin{equation}\label{eq:quo_SUM}
	\dot x_{[i]}= A x_{[i]}+\sum_k w_{[i],k}^- \rep (p_{[j]})	\left(
	\begin{array}{l}
	1 \\
	0
	\end{array}
	\right)
	+\sum_k w_{[i],k}^+ \act (p_{[j]})	\left(
	\begin{array}{l}
	1 \\
	0
	\end{array}
	\right), 
	\end{equation}
	where $[j]:=C_k$ if $j\in C_k$ and $k$ is such that $I_i^\pm\cap C_k\ne \emptyset$. Notice that by  (\ref{cond:synSUM+})-(\ref{cond:synSUM-}), $I_i^\pm\cap C_k\ne \emptyset$ if and only if $I_{l}^\pm\cap C_k\ne \emptyset$ for any $l\in [i]$.	 It follows that  $\Delta_P$ is flow-invariant and consequently, a synchrony space for the SUM equations model (\ref{eq:SUM}).
	
	Assume that  $\Delta_P$ is a synchrony space for (\ref{eq:SUM}). Then, $\Delta_P$ is flow-invariant. In particular, we have $I_i^\pm\cap C_k\ne \emptyset$ if and only if $I_{l}^\pm\cap C_k\ne \emptyset$ for any $i,l$ in the same part of P.  Moreover, the flow restricted to $\Delta_P$ is of form (\ref{eq:quo_SUM}) and for $i,l$ in the same part of P, we have
	\[\sum_{k=1}^m w_{i,k}^- \rep (p_{[j]})
	+\sum_{k=1}^m w_{i,k}^+ \act (p_{[j]})=
	\sum_{k=1}^m w_{l,k}^- \rep (p_{[j]})
	+\sum_{k=1}^m w_{l,k}^+ \act (p_{[j]}).\]
For $k\in \{1,\dots, m\}$, let $p_{[j]}=0$ for all $[j]\ne C_{k}$.
For $[j]=C_k$, if we let $p_{[j]}\to 0$, then $\rep  (p_{[j]})\to 1$ and  $ \act  (p_{[j]})\to 0$ hold, which implies that $w_{i,k}^-=w_{l,k}^-$ 
(cf. (\ref{eq;reg_Hill_like_0})). Similarly, if we let  $p_{[j]}\to \infty$ for $[j]=C_k$, then  $w_{i,k}^+=w_{l,k}^+$  
(cf. (\ref{eq;reg_Hill_like_infty})).  
\end{proof}

	\begin{rem}\rm
		It follows from Proposition~\ref{prop:synSUM} that, given a partition $P$ of the gene set of a GRN, the associated polydiagonal susbspace $\Delta_P$ is a synchrony subspace for the SUM model of the GRN if and only if it is left invariant by the weighted adjacency matrices $W^+$ and $W^-$ of the GRN. We have then, using the work of Aguiar and Dias in \cite{AD14} and \cite{AD18}, that the set of the synchrony subspaces (synchrony patterns) for the SUM model of a GRN can be computed using  the algorithm in Section 6 of \cite{AD14}. The algorithm can be executed to find the set of polydiagonal subspaces that are left invariant by one of the adjacency matrices $W^+$ or $W^-$ and then the set of synchrony subspaces of the GRN is the subset of those polydiagonals that are also left invariant by the other adjacency matrix.
\hfill $\Diamond$	
	\end{rem}

\begin{prop}\label{prop:SUMquotyes}
Take an $n$-gene GRN with adjacency matrices $W^+$ and $W^{-}$. Let $P$ be a synchronization partition of the gene set with classes $C_1, \ldots, C_m$ refining the gene equivalence partition. Consider the $m \times m$ matrices $Q^{+} = \left[ q^{+}_{ik} \right]$ and   $Q^{-} = \left[ q^{-}_{ik} \right]$ where
$$ 
\sum_{j \in I_i^+ \cap C_k} w^+_{ij} =  q^{+}_{ik}  \quad \mbox{ and } 
\sum_{j \in I_i^{-} \cap C_k} w^{-}_{ij} = q^{-}_{ik}\quad \left(i,k= 1, \dots, m\right)\, . 
$$
 The restriction of SUM equations model (\ref{eq:SUM}) to $\Delta_P$ is a SUM equations model consistent with the 
quotient GRN of $m$ genes with activation and repression weighted adjacency matrices $Q^+$ and $Q^{-}$, respectively.
\end{prop}

\begin{proof}
Follows trivially from the conditions in Proposition~\ref{prop:synSUM} and its proof. 
\end{proof}

\begin{figure}
	\begin{center}
\begin{tikzpicture}
 [scale=.15,auto=left, node distance=1.5cm, 
 ]
\node[fill=white,style={circle,draw}] (n1) at (24,10) {\small{1}};
\node[fill=white,style={circle,draw}] (n2) at (24,25) {\small{2}};
\node[fill=white,style={circle,draw}] (n3) at (4,10) {\small{3}};
\node[fill=white,style={circle,draw}] (n4) at (4,25) {\small{4}};
\node[fill=white,style={circle,draw}] (n5) at (14,0) {\small{5}};

\node[fill=white,style={circle,draw}] (n6) at (44,15) {\small{[1]}};
\node[fill=white,style={circle,draw}] (n7) at (54,15) {\small{[4]}};

 \path 
 (n1) [->] edge[bend left=10,thick] node [near end, right] {{\tiny $2$}} (n2)
 (n1) [-|] edge[bend left=10,thick] node [near end, below=0.1pt] {{\tiny $3$}} (n5)
 (n2) [->] edge[bend left=20,thick] node [near end, right] {{\tiny $2$}} (n1)
 (n2) [->] edge[thick] node [near end, above] {{\tiny $1$}} (n3)
(n2) [-|] edge[bend left=10,thick] node [near end, below]  {{\tiny $2$}} (n4)
 (n3) [-|] edge[thick] node [ above=0.1pt] {{\tiny $1$}} (n1)
  (n3) [-|] edge[thick, loop left] node  {{\tiny $1$}} (n3)
  (n3) [-|] edge[thick] node [near end, left] {{\tiny $1$}} (n4)
 (n4) [-|] edge[thick] node [near end, above=0.1pt] {{\tiny $4$}} (n1)
 (n4) [-|] edge[bend left=10,thick] node [near end, above] {{\tiny $5.5$}} (n2)
  (n5) [-|] edge[bend left=10,thick] node [near end, above] {{\tiny $5$}} (n1)
   (n5) [-|] edge[thick] node [near end, above] {{\tiny $3.5$}} (n2)
    (n5) [-|] edge[thick] node [near end, below=0.1pt] {{\tiny $9$}} (n3)
    
      (n6) [->] edge[thick, loop left] node  {{\tiny $2$}} (n6)
      (n6) [-|] edge[bend left=20,thick] node [near end, above] {{\tiny $3$}} (n7)
           (n7) [-|] edge[bend left=20,thick] node [near end, below] {{\tiny $9$}} (n6)
   ;
\end{tikzpicture} 
\caption{The GRNs of Example~\ref{ex:5x5}: (left) a five gene GRN and (right) its quotient network corresponding to the partition $P =\left\{ [1]=\{1,2,3\},\, [4]=\{ 4,5\}\right\}$.}
\label{fig:ex5x5}
\end{center}
\end{figure}
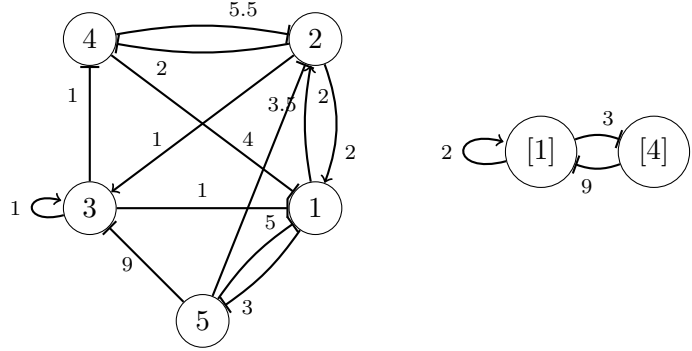

\begin{exam} \label{ex:5x5}
Take the five gene GRN on the left of Figure~\ref{fig:ex5x5} with the following $5 \times 5$ activation and repression weighted adjacency matrices:
$$
W^+ = 
\left(
\begin{array}{rrr|rr}
0& 2 & 0 & 0 & 0\\
2 & 0 & 0 & 0 & 0 \\
0 & 1 & 1 & 0 & 0 \\
\hline
0 & 0 & 0 & 0 & 0 \\
0 & 0 & 0 & 0 & 0
\end{array}
\right), \quad 
W^{-} = 
\left(
\begin{array}{rrr|rr}
0& 0 & 0 & 4 & 5\\
0 & 0 & 0 & 5.5 & 3.5 \\
0 & 0 & 0 & 0 & 9 \\
\hline
0 & 2 & 1 & 0 & 0 \\
3 & 0 & 0 & 0 & 0
\end{array}
\right)\, .
$$

By Proposition~\ref{prop:synSUM}, the partition $P =\left\{ [1]=\{1,2,3\},\, [4]=\{ 4,5\}\right\}$ is a synchronization partition and corresponds to the synchrony space $\Delta_P = \{ x:\, x_1 = x_2 = x_3,\, x_4 = x_5\}$ for the SUM model equations if $P$ refines the gene equivalence partition. By Proposition~\ref{prop:SUMquotyes}, 
the corresponding quotient network is the two-gene GRN on the right of Figure~\ref{fig:ex5x5} with the following  activation and repression weighted adjacency matrices: 
\begin{equation} \label{eq:exquots}
Q^+ = 
\left(
\begin{array}{r|r}
2 & 0 \\
\hline
0 & 0 
\end{array}
\right), \quad 
Q^{-} = 
\left(
\begin{array}{r|r}
0 & 9\\
\hline
3 & 0
\end{array}
\right)\, .
\end{equation}
Note that $\Delta_P$ is not a synchrony space for the MULT model equations. 
\hfill $\Diamond$
\end{exam}

\begin{rem}
For any five gene GRN with $5 \times 5$ activation and repression weighted adjacency matrices given by 
$$
W^+ = 
\left(
\begin{array}{rrr|rr}
0& 2 & 0 & 0 & 0\\
2 & 0 & 0 & 0 & 0 \\
0 & 1 & 1 & 0 & 0 \\
\hline
0 & 0 & 0 & 0 & 0 \\
0 & 0 & 0 & 0 & 0
\end{array}
\right), \quad 
W^{-} = 
\left(
\begin{array}{rrr|rr}
0& 0 & 0 & w^{-}_{14} & w^{-}_{15} \\
0 & 0 & 0 & w^{-}_{24} & w^{-}_{25} \\
0 & 0 & 0 & 0 & 9 \\
\hline
0 & 2 & 1 & 0 & 0 \\
3 & 0 & 0 & 0 & 0
\end{array}
\right)
$$
where the weights satisfy $w^{-}_{14} + w^{-}_{15} = w^{-}_{24} + w^{-}_{25} = 9$, if the gene set partition $P$ with classes $ \{1,2,3\},\, \{ 4,5\}$ refines the gene equivalence partition, then $P$ is a synchronization partition corresponding  to the synchrony space $\Delta_P = \{ x:\, x_1 = x_2 = x_3,\, x_4 = x_5\}$  for the SUM model equations. The corresponding quotient network of two-gene GRN has activation and repression weighted adjacency matrices given in (\ref{eq:exquots}). This follows from Propositions~\ref{prop:synSUM} and \ref{prop:SUMquotyes}. This example illustrates that there is an infinite number of five gene GRNs admitting the synchronization partition $P$ and leading to the same two-gene quotient GRN for the SUM models. 
\hfill $\Diamond$
\end{rem}

\begin{rem}\rm  \label{rmk:funny} 
Under the assumptions of Proposition~\ref{prop:synSUM} and recall equation (\ref{eq:quo_SUM}),
if $w^{+}_{i,k} = w^{-}_{i,k}$, then as $ \rep (p_k) = 1- \act (p_k)$, the influence of  the $k$-th cluster $C_k$ to the evolution of the $i$-th gene reduces to $w^{+}_{i,k} \act (p_k) + w^{-}_{i,k}  \rep (p_k) =  w^{-}_{i,k} = w^{+}_{i,k}$,  which becomes independent on the variable $p_k$. See Example \ref{ex:funny}.   
For general weighted networks,  Aguiar, Dias and Ferreira~\cite{ADF17} point out a similar phenomenon and in that case $\Delta_P$ is called a {\it spurious} synchrony pattern. See Definition 2.9 of \cite{ADF17}.
\hfill $\Diamond$
\end{rem}

\begin{rem}
There are GRNs such that  $I_i^+ \cap I_i^- = \emptyset$, for all genes $i$, and admitting a synchronization partition such that, in the quotient, the same gene activates and represses another gene, as the following example illustrates.
\hfill $\Diamond$
\end{rem}

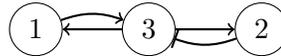
\begin{figure}[h!]
\begin{center}
\begin{tikzpicture}
 [scale=.15,auto=left, node distance=1.5cm, 
 ]
\node[fill=white,style={circle,draw}] (n1) at (4,0) {\small{1}};
\node[fill=white,style={circle,draw}] (n2) at (24,0) {\small{2}};
\node[fill=white,style={circle,draw}] (n3) at (14,0) {\small{3}};
 \path 
 (n1) [->] edge[thick, bend left=20] node [near end, above] {} (n3)
 (n2) [-|] edge[thick, bend left=20] node [near end, above] {} (n3)
 (n3) [->] edge[thick] node [near end, above] {} (n1)
 (n3) [->] edge[thick] node [near end, above] {} (n2)
   ;
		\end{tikzpicture}	
		\end{center}
\caption{GRN with three genes.}
\label{fig:funny}
\end{figure}

\begin{figure}[h!]
\begin{center}
\begin{tikzpicture}
 [scale=.15,auto=left, node distance=1.5cm, 
 ]
\node[fill=white,style={circle,draw}] (n1) at (4,0) {\small{[1]}};
\node[fill=white,style={circle,draw}] (n3) at (14,0) {\small{[3]}};
 \path 
 (n1) [->] edge[thick, bend left=20] node [near end, above] {} (n3)
 (n1) [-|] edge[thick, bend right=20] node [near end, above] {} (n3)
 (n3) [->] edge[thick] node [near end, above] {} (n1)
   ;
		\end{tikzpicture}	
		\end{center}
\caption{GRN with two genes where gene $[1]$ activates and supresses gene $[3]$ with the same strenght.}
\label{fig:qfunny}
\end{figure}
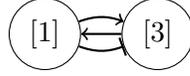

\begin{exam}\label{ex:funny}
Consider the three-gene GRN in Figure~\ref{fig:funny} with $3 \times 3$ activation and repression weighted adjacency matrices given by 
$$
W^+ = 
\left(
\begin{array}{rr|r}
0 & 0 & 1\\
 0 & 0 & 1 \\
 \hline 
  1 & 0 & 0 
\end{array}
\right), \quad 
W^{-} = 
\left(
\begin{array}{rr|r}
0& 0 & 0 \\
0 & 0 & 0 \\
\hline 
0 & 1 & 0
\end{array}
\right)\, .
$$
By Proposition~\ref{prop:synSUM}, assuming genes $1$ and $2$ are equivalent, we have that  $P=\{ C_1 = \{1,2\}, C_2 = \{3\}\}$ is a synchronization partition for the SUM equations models. 
Now, the two-gene quotient GRN as defined in Proposition~\ref{prop:SUMquotyes}, with 
 $2 \times 2$ activation and repression weighted adjacency matrices given by 
$$
Q^+ = 
\left(
\begin{array}{r|r}
0 & 1 \\
\hline
1 & 0 
\end{array}
\right), \quad 
Q^{-} = 
\left(
\begin{array}{r|r}
0 & 0\\
\hline
1 & 0
\end{array}
\right),$$
 and satisfies $I_{[3]}^+ \cap I_{[3]}^- = \{ [1]\} \ne \emptyset$. That is, in the quotient, gene $[1]$ activates and represses gene $[3]$, see Figure~\ref{fig:qfunny}. This follows from the fact that $I_3^+ = \{1\}$, $I_3^{-} = \{ 2\}$ and $1,2$ form the part $C_1$ of the synchronization partition $P$.   Equations for the SUM model  
 $$
 \left\{
\begin{array}{l}
\dot{x}_1 = 
A_1 x_1 + \act (p_3) 
\left(
\begin{array}{l}
1 \\
0
\end{array}
\right)
\vspace{2mm} 
\\
\dot{x}_2 = 
A_1 x_2 + \act (p_3) 
\left(
\begin{array}{l}
1 \\
0
\end{array}
\right)
\\
\dot{x}_3 = 
A_3 x_3 + \left[ 
 \act (p_1) + \rep (p_2) \right] 
\left(
\begin{array}{l}
1 \\
0
\end{array}
\right)
\end{array}
\right. 
 $$
restricted to $\Delta_{P} = \{ x:\, x_1 = x_2\}$ are given by:
 $$
 \left\{
\begin{array}{l}
\dot{x}_1 = 
A_1 x_1 + \act (p_3) 
\left(
\begin{array}{l}
1 \\
0
\end{array}
\right)
\\
\dot{x}_3 = 
A_3 x_3 + 
\left[ \act (p_1) + \rep (p_1) \right] 
\left(
\begin{array}{l}
1 \\
0
\end{array}
\right)
\end{array}
\right. \, . 
$$
With the assumption $\act (p_1) + \rep (p_1) \equiv 1$, these equations simplify to
 $$
  \left\{
\begin{array}{l}
\dot{x}_1 = 
A_1 x_1 + \act (p_3) 
\left(
\begin{array}{l}
1 \\
0
\end{array}
\right)
\\
\dot{x}_3 = 
A_3 x_3 + 
\left(
\begin{array}{l}
1 \\
0
\end{array}
\right)
\end{array}
\right. \, .
$$
Thus, the gene equation for $[3]$ becomes independent of gene $[1]$. This example also illustrates Remark~\ref{rmk:funny}. 
\hfill $\Diamond$
\end{exam}

\subsection{Synchrony for the MULT model}

Let $W^\pm=(w_{ij}^\pm)$ and $M^\pm=(m_{ij}^\pm)$  be, respectively, the weighted matrices and the multiplicity matrices of the MULT model (\ref{eq:MULT}).

\begin{prop}\label{prop:MULTsyn}
	Take an $n$-gene GRN with adjacency matrices $W^\pm$ and $M^\pm$ and a partition $P$ of the gene set into classes $C_1, \ldots, C_m$ refining the gene equivalence class. The partition $P$ corresponds to a synchrony pattern for the MULT equations model (\ref{eq:MULT}) if and only if for each part $C$ of $P$,  we have that: \\
	(i) for $k=1, \ldots, m$, 
	$$\sum_{j\in I_i^-\cap C_k} m_{ij}^-\quad  \mbox{ is constant for } i \in C$$
	and 
	$$\sum_{j\in I_i^+\cap C_k} m_{ij}^+\quad  \mbox{ is constant for } i \in C$$
	(ii) \begin{equation} \label{eq:prodcondition} \prod_{j \in I^+_i} w^+_{ij} \prod_{j \in I^{-}_i} w^{-}_{ij} \mbox{ is constant for $i \in C$.} \end{equation}  
\end{prop}
\begin{proof} 
Assume that (i)-(ii) hold on any equivalence class $C$ under partition $P$. Denote by $\mathbf{m}_{i,k}^\pm$ the  sum of the multiplicities $m_{ij}^\pm$ in $I_i^\pm\cap C_k$, that is, $\mathbf{m}_{i,k}^{-} = \sum_{j\in I_i^-\cap C_k} m_{ij}^-$ and $\mathbf{m}_{i,k}^{+} = \sum_{j\in I_i^+\cap C_k} m_{ij}^+$. Then, we have for any $i\in C$ that
\begin{align*}
\dot x_i&=Ax_i+
\prod_{j \in I^{-}_i} w^{-}_{ij} \rep^{m_{ij}^-} (p_j) \, 
\prod_{j \in I^{+}_i} w^{+}_{ij} \act^{m_{ij}^+} (p_j) 
\left(
\begin{array}{l}
1 \\
0
\end{array}
\right) \\
&=Ax_i+
\big(\underbrace{\prod_{j \in I^{-}_i} w^{-}_{ij}\prod_{j \in I^{+}_i} w^{+}_{ij}}_{w_i}\big) \prod_k\prod_{j\in I_i^-\cap C_k} \rep^{m_{ij}^-} (p_{j}) \, 
 \prod_k\prod_{j\in I_i^+\cap C_k} \act^{m_{ij}^+}  (p_{j}) 
\left(
\begin{array}{l}
1 \\
0
\end{array}
\right),
\end{align*}	
where the product $w_i$, by (ii), remains the same in the equivalence class of $i$. In fact, $w_i$ is the product of the $w_{ij}^{\pm}$ where $j$ runs through $I_i^\pm$.
Moreover, by (i), the sums $\mathbf{m}_{i,k}^{\pm}$ are also constant in the equivalence class of $i$.
Thus, for any $i\in C$, the dynamics of $x_{[i]}\in \Delta_P$ is governed by
\begin{equation}\label{eq:MULT_quot}
\dot x_{[i]}=Ax_{[i]}+w_{[i]} \prod_k \rep^{\mathbf{m}_{[i],k}^-} (p_{[j]}) \, 
\prod_k \act^{\mathbf{m}_{[i],k}^+} (p_{[j]}) \left(
\begin{array}{l}
1 \\
0
\end{array}
\right),
\end{equation}
where $[j]:=C_k$ for $j\in C_k$ and $k$ is such that $I^{\pm}_i\cap C_k\ne  \emptyset$. 
By (ii), $I^{\pm}_i\cap C_k\ne  \emptyset$ if and only if $I^{\pm}_l\cap C_k\ne  \emptyset$ for any $l\in[i]$. Thus, $\Delta_P$ is flow-invariant and a synchrony space for the MULT equations model  (\ref{eq:MULT}).

Assume that $\Delta_P$ is a synchrony space for  (\ref{eq:MULT}). Then, $\Delta_P$ is flow-invariant. In particular, $I_i^\pm\cap C_k\ne \emptyset$ if and only if $I_l^\pm\cap C_k\ne \emptyset$ for any $i,l$ in the same part of $P$. Moreover, the flow restricted to $\Delta_P$ is of shape (\ref{eq:MULT_quot}) and for any $i,l$ in the same part of $P$, we have
\begin{equation}\label{eq:pf_mult_2}
w_i \prod_k \rep^{\mathbf{m}_{[i],k}^-} (p_{[j]}) \, 
\prod_k \act^{\mathbf{m}_{[i],k}^+} (p_{[j]})
=w_l \prod_k \rep^{\mathbf{m}_{[l],k}^-} (p_{[j]}) \, 
\prod_k \act^{\mathbf{m}_{[l],k}^+} (p_{[j]}).
\end{equation}
Assume to the contrary of (i) that $\mathbf{m}_{i,k}^-\ne \mathbf{m}_{l,k}^-$ for some $k$. Then, by letting $p_{[j]}=0$ for all $[j]$ except when $[j]=C_k$ in (\ref{eq:pf_mult_2}), we have ${\mathrm rep}^M(p){\mathrm act}^N(p)\equiv \text{constant}$ for $p=p_{[j]}$ and some positive integers $M,N$. But since ${\mathrm rep}^M(p){\mathrm act}^N(p)\to 0$ both when $p\to 0$ and $p\to \infty$  by (\ref{eq;reg_Hill_like_0})-(\ref{eq;reg_Hill_like_infty}), we must have  ${\mathrm rep}^M(p){\mathrm act}^N(p)\equiv 0$, which contradicts to the fact that ${\mathrm rep}$ and ${\mathrm act}$ are non-zero functions. Therefore, $\mathbf{m}_{i,k}^\pm=\mathbf{m}_{l,k}^\pm$ and (i) follows. By (\ref{eq:pf_mult_2}), we conclude that $w_i=w_l$ and (ii) follows.
\end{proof}

	\begin{rem}\rm
		Given a GRN with weighted adjacency matrices $W^+$ and $W^-$, consider the diagonal matrix
		$W^*$ where the $ii$ entry is given by the product of the non-zero entries of the $i$-th row of both matrices $W^+$ and $W^-$,or zero otherwise.
		From Proposition~\ref{prop:MULTsyn}, given a partition $P$ of the gene set of the GRN, the associated polydiagonal susbspace $\Delta_P$ is a synchrony subspace for the MULT model of the GRN if and only if it is left invariant by the multiplicity matrices $M^+$ and $M^-$ of the GRN and by the matrix $W^*$. We have then, using the work of Aguiar and Dias in \cite{AD14} and \cite{AD18}, that the set of the synchrony patterns for the MULT model of a GRN can be computed using  the algorithm in Section 6 of \cite{AD14}. The algorithm can be executed to find the set of polydiagonal subspaces that are left invariant by one of the multiplicity matrices $M^+$ or $M^-$ and then the set of synchrony subspaces of the GRN is the subset of those polydiagonals that are also left invariant by the other multiplicity matrix and by the diagonal matrix $W^*$.
		\hfill $\Diamond$ 
	\end{rem}

\begin{prop}\label{prop:MULT_quo}
	Take an $n$-gene GRN with multiplicity matrices $M^\pm$ and weight matrices $W^\pm$ and a synchronized partition $P$ of the gene set into classes $C_1, \ldots, C_m$ refining the gene equivalence class for MULT equations model (\ref{eq:MULT0}). The restriction of (\ref{eq:MULT0}) to the synchrony space $\Delta_P$ is a MULT equations model consistent with any  quotient $m$-gene GRN with the $m \times m$ multiplicity matrices $N^+=\left[ n^{+}_{ik} \right]$ and $N^-=\left[ n^{-}_{ik} \right]$  given by
$$n_{ik}^+=\sum_{j\in I_i^+\cap C_k} m_{ij}^+,\quad 	n_{ik}^-=\sum_{j\in I_i^-\cap C_k} m_{ij}^-, \quad \left(i,k= 1, \dots, m\right) $$
and $m \times m$ activation and repression weighted adjacency matrices $Q^{+} = \left[ q^{+}_{ik} \right]$ and   $Q^{-} = \left[ q^{-}_{ik} \right]$ satisfying 
\begin{equation}	\label{eq:prodcond}
\prod_{k=1}^m q^{+}_{ik} \, q^{-}_{ik}  = \prod_{k=1}^m \left(\prod_{j \in I_i^+ \cap C_k} w^+_{ij}\right) \left( \prod_{j \in I_i^{-} \cap C_k} w^{-}_{ij}\right)  \quad \left(i= 1, \dots, m\right)\, . \end{equation}
Here, each product $\prod_{j \in I_i^+ \cap C_k} w^+_{ij}$ is considered only for $k$ such that  $I_i^+ \cap C_k \not=\emptyset$. In that case, $\prod_{j \in I_i^+ \cap C_k} w^+_{ij}$ is positive. 
Similarly, each product $\prod_{j \in I_i^{-} \cap C_k} w^{-}_{ij}$ is taken for $k$ such that $I_i^{-} \cap C_k \not=\emptyset$.  
\end{prop}

\begin{proof}
	This follows from Proposition~\ref{prop:MULTsyn} and the definition of MULT equations model (\ref{eq:MULT0}).
\end{proof}

\begin{rem}
Note that, in general, there is no uniqueness on the \textcolor{red}{quotient} GRN associated with the MULT model equations restricted to a synchrony space. 
This  lack of uniqueness is due to the dependence of the product conditions (\ref{eq:prodcond}), defining the activation and repression weighted adjacency matrices. See Example \ref{ex:quo_not_prod} below. 
{Thus, the choice of a particular quotient GRN can be made taking into account the specifics of the problem under analysis.}
\hfill $\Diamond$
\end{rem}

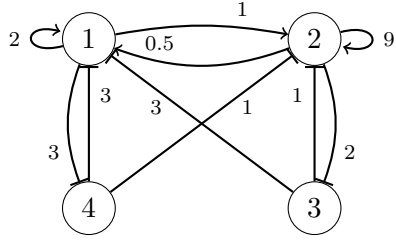
\begin{figure}
	\begin{center}
\begin{tikzpicture}
 [scale=.15,auto=left, node distance=1.5cm, 
 ]
\node[fill=white,style={circle,draw}] (n1) at (4,25) {\small{1}};
\node[fill=white,style={circle,draw}] (n2) at (24,25) {\small{2}};
\node[fill=white,style={circle,draw}] (n3) at (24,10) {\small{3}};
\node[fill=white,style={circle,draw}] (n4) at (4,10) {\small{4}};

 \path 
   (n1) [->] edge[thick, loop left] node  {{\tiny $2$}} (n1)
 (n1) [->] edge[bend left=10,thick] node [near end, above] {{\tiny $1$}} (n2)
  (n1) [-|] edge[bend right=20,thick] node [near end, left] {{\tiny $3$}} (n4)
   
 (n2) [->] edge[thick, loop right] node  {{\tiny $9$}} (n2)
 (n2) [->] edge[bend left=20,thick] node [near end, above] {{\tiny $0.5$}} (n1)
 (n2) [-|] edge[bend left=20,thick] node [near end, right] {{\tiny $2$}} (n3)

 (n3) [-|] edge[thick] node [near end,below] {{\tiny $3$}} (n1)
(n3) [-|] edge[thick] node [near end, left] {{\tiny $1$}} (n2)
  
 (n4) [-|] edge[thick] node [near end, right] {{\tiny $3$}} (n1)
 (n4) [-|] edge[thick] node [near end, below] {{\tiny $1$}} (n2)
 
   ;
\end{tikzpicture} 
\caption{The four gene GRN of Example~\ref{ex:quo_not_prod}.}
\label{fig:ex4x4}
\end{center}
\end{figure}

\begin{exam}\label{ex:quo_not_prod}
Take the 4-gene GRN in Figure~\ref{fig:ex4x4} with the following $4 \times 4$ activation and repression weighted adjacency matrices, 
$$
W^+ = 
\left(
\begin{array}{rr|r|r}
2& 0.5 & 0 &  0\\
1 & 9 & 0 & 0 \\
\hline
0 & 0 & 0 & 0 \\
\hline 
0 & 0 & 0 & 0 
\end{array}
\right), \  
W^{-} = 
\left(
\begin{array}{rr|r|r}
0& 0 & 3 & 3\\
0 & 0 & 1 & 1 \\
\hline
0 & 2 & 0 & 0 \\
\hline 
3 & 0 & 0 & 0
\end{array}
\right),$$ 
and multiplicity matrices,
$$M^+ = 
\left(
\begin{array}{rr|r|r}
1& 1 & 0 &  0\\
1 & 1 & 0 & 0 \\
\hline
0 & 0 & 0 & 0 \\
\hline 
0 & 0 & 0 & 0 
\end{array}
\right), \  
M^{-} = 
\left(
\begin{array}{rr|r|r}
0& 0 & 1 & 1\\
0 & 0 & 1 & 1 \\
\hline
0 & 1 & 0 & 0 \\
\hline 
1 & 0 & 0 & 0
\end{array}
\right)\, .
$$
Assume that genes $1,2$ have the same internal dynamics. Consider the gene set partition $P =\left\{ C_1 = \{1,2\},\, C_2 = \{3\},\, C_3 = \{ 4\}\right\}$.  By Proposition~\ref{prop:MULTsyn}, we have that  $P$ is a synchronization partition. Note that $2 \times 0.5\times  3\times  3 = 1\times 9\times 1\times 1 =9$ and so condition (\ref{eq:prodcondition}) is satisfied. The synchronization partition $P$ corresponds to the synchrony space $\Delta_P = \{ x:\, x_1 = x_2\}$ for the MULT (and PROD) model equations:
\begin{equation}\label{eq:yesPROD}
\left\{
\begin{array}{l}
\dot{x}_1 = 
A_1 x_1 + \left[ 2 \act (p_1)  \, 0.5  \act (p_2)  \, 3 \rep (p_3)  \, 3  \rep (p_4) \right]
\left(
\begin{array}{l}
1 \\
0
\end{array}
\right)
\\
\dot{x}_2 = 
A_1 x_2 + \left[  \act (p_1)  \, 9  \act (p_2)  \,  \rep (p_3)  \,   \rep (p_4) \right]
\left(
\begin{array}{l}
1 \\
0
\end{array}
\right)
\\
\dot{x}_3 = 
A_3 x_3 + 
 2 \rep (p_2) 
\left(
\begin{array}{l}
1 \\
0
\end{array}
\right)\\
\dot{x}_4 = 
A_4 x_4 + 
 3 \rep (p_1) 
\left(
\begin{array}{l}
1 \\
0
\end{array}
\right)
\end{array}
\right.
\end{equation}

The restriction of  equations (\ref{eq:yesPROD}) to $\Delta_P = \{ x:\, x_1 = x_2\}$  is:
\begin{equation}\label{eq:notPROD}
\left\{
\begin{array}{l}
\dot{x}_1 = 
A_1 x_1 + \left[ 9 \left( \act (p_1)\right)^2 \,   \rep (p_3)  \, \rep (p_4) \right]
\left(
\begin{array}{l}
1 \\
0
\end{array}
\right)
\\
\dot{x}_3 = 
A_3 x_3 + 
 2  \rep (p_1) 
\left(
\begin{array}{l}
1 \\
0
\end{array}
\right)\\
\dot{x}_4 = 
A_4 x_4 + 
 3 \rep (p_1) 
\left(
\begin{array}{l}
1 \\
0
\end{array}
\right)
\end{array}
\right. \, .
\end{equation}
 These equations are MULT model equations for 
   any 3-gene GRN with the following $3 \times 3$ multiplicity matrices,
$$N^+ = 
\left(
\begin{array}{r|r|r}
2 & 0 &  0\\
\hline
0 & 0 & 0 \\
\hline 
0 & 0 & 0 
\end{array}
\right), \  
N^{-} = 
\left(
\begin{array}{r|r|r}
0 & 1 & 1\\
\hline
1 & 0 & 0 \\
\hline 
1 & 0 & 0
\end{array}
\right) 
$$
	and activation and repression weighted adjacency matrices,
$$
Q^+ = 
\left(
\begin{array}{r|r|r}
q^+_{11}& 0 &  0\\
\hline
0  & 0 & 0 \\
\hline 
0  & 0 & 0 
\end{array}
\right), \  
Q^{-} = 
\left(
\begin{array}{r|r|r}
 0 & q^-_{12} & q^-_{13} \\
\hline
2 & 0 & 0 \\
\hline 
3 & 0 & 0
\end{array}
\right),$$
such that $q^+_{11} q^-_{12} q^-_{13} =9$.
Two particular choices of activation and repression weighted adjacency matrices are given by, respectively, 
	$$
	Q^+ = 
	\left(
	\begin{array}{r|r|r}
	1& 0 &  0\\
	\hline
	0  & 0 & 0 \\
	\hline 
	0  & 0 & 0 
	\end{array}
	\right), \  
	Q^{-} = 
	\left(
	\begin{array}{r|r|r}
	0 & 3 & 3\\
	\hline
	2 & 0 & 0 \\
	\hline 
	3 & 0 & 0
	\end{array}
	\right).$$
and
	$$
	Q^+ = 
	\left(
	\begin{array}{r|r|r}
	9& 0 &  0\\
	\hline
	0  & 0 & 0 \\
	\hline 
	0  & 0 & 0 
	\end{array}
	\right), \  
	Q^{-} = 
	\left(
	\begin{array}{r|r|r}
	0 & 1 & 1\\
	\hline
	2 & 0 & 0 \\
	\hline 
	3 & 0 & 0
	\end{array}
	\right).$$

	Note that equations (\ref{eq:notPROD}) fit the MULT model but do not fit the PROD model. 
\hfill $\Diamond$
\end{exam}

As it is pointed out by Example \ref{ex:quo_not_prod} above, the PROD model shows some weaknesses in the context of synchrony spaces and their quotient networks, which are indispensable in describing robust synchrony patterns accommodated by GRN network structure. This disadvantage ceases to exist in the more general MULT model.

 \section{Network lifting enumeration and gene duplication} \label{sec:enumlift}
  
 As stated in the introductory Section~\ref{sec:intro}, for general coupled cell networks and weighted networks,  fixing a $k$-cell (quotient) network, if $n>k$ then there are many $n$-cell networks which are lifts (or inflations) of the $k$-cell network, that is, that admit the fixed network as a quotient. Moreover, there is a method of enumerating the lifts of a fixed $k$-cell network. See Theorem 2.5 of Aguiar, Dias, Golubitsky and Leite~\cite{ADGL09} valid for coupled cell systems following the formalisms of Stewart, Golubitsky and Pivato~\cite{SGP03}, Golubitsky, Stewart and T\"{o}r\"{o}k~\cite{GST05} and Field~\cite{F04}, and its extension to weighted networks for coupled cell systems with additive input structure in  Theorem 2.13 of Aguiar and Dias~\cite{AD18}. The enumeration method relies upon the characterization of the network adjacency matrices of the larger networks determined by the smaller network adjacency matrices. Also, fixing $n$, the number of $n$-cell lifts is finite for nonnegative integer matrices and nonfinite in the weighted setup.

 In this section, we combine Theorem 2.13 of Aguiar and Dias~\cite{AD18} with Propositions~\ref{prop:synSUM}, \ref{prop:SUMquotyes}, and 
 Propositions~\ref{prop:MULTsyn}, \ref{prop:MULT_quo}, obtaining a characterization method of the $n$-gene GRNs that are lifts of a fixed $k$-gene GRN, for both the SUM and MULT models, using the activation and regression weighed matrices (and the multiplicities matrices for the MULT model quotients).

 \begin{thm}\label{thm:enumSUM}
 Consider $Q$ an $m$-gene GRN for the SUM equations model with activation and repression weighted adjacency matrices $Q^{+} = \left[ q^{+}_{ik} \right]$ and   $Q^{-} = \left[ q^{-}_{ik} \right]$. An $n$-gene GRN for the SUM equations model with set of genes $\{1, \ldots, n\}$, where $n > m$, is a lift of $Q$ if and only if there is a partition of $\{1, \ldots, n\}$ into $m$ classes, $C_1, \ldots, C_m$, refining the gene equivalence partition  such that, after renumbering the genes if necessary, the activation and repression weighted adjacency matrices $W^+$ and $W^{-}$ have the following block structures: 
 \begin{equation}\label{eq:enumSUM}
 W^{+} = 
 \left( 
 \begin{array}{ccc}
 Q_{11}^{+}  & \cdots & Q_{1m}^{+}  \\
 \vdots & \cdots & \vdots \\
 Q_{m1}^{+}  & \cdots & Q_{mm}^{+} 
 \end{array}
 \right) ,\qquad 
 W^{-} = 
 \left( 
 \begin{array}{ccc}
 Q_{11}^{-}  & \cdots & Q_{1m}^{-}  \\
 \vdots & \cdots & \vdots \\
 Q_{m1}^{-}  & \cdots & Q_{mm}^{-} 
 \end{array}
 \right)
 \end{equation}
 where each $Q_{ik}^{\pm}$ is an $\# C_i \times \# C_k$-matrix with nonnegative real entries whose  row  sum is $q^{\pm}_{ik}$.   
 \end{thm}

\begin{proof}
Direct application of Proposition~\ref{prop:synSUM} and Theorem 2.13 in \cite{AD18}. 
\end{proof}


\begin{thm}\label{thm:enumMULT} 
 Consider $Q$ an $m$-gene GRN for the MULT equations model with $m \times m$ multiplicity matrices $N^+=\left[ n^{+}_{ik} \right]$, $N^-=\left[ n^{-}_{ik} \right]$ and activation and repression weighted adjacency matrices $Q^{+} = \left[ q^{+}_{ik} \right]$ and   $Q^{-} = \left[ q^{-}_{ik} \right]$.  An $n$-gene GRN for the MULT equations model with set of genes $\{1, \ldots, n\}$, where $n > m$, is a lift of $Q$ if and only if there is a partition of $\{1, \ldots, n\}$ into $m$ classes, $C_1, \ldots, C_m$, refining the gene equivalence partition  such that, after renumbering the genes if necessary: \\
 (i) The $n \times n$ multiplicity matrices $M^+$ and $M^{-}$ have the following block structures: 
 \begin{equation}\label{eq:enumMULT}
 M^{+} = 
 \left( 
 \begin{array}{ccc}
 N_{11}^{+}  & \cdots & N_{1m}^{+}  \\
 \vdots & \cdots & \vdots \\
 N_{m1}^{+}  & \cdots & N_{mm}^{+} 
 \end{array}
 \right) ,\qquad 
 M^{-} = 
 \left( 
 \begin{array}{ccc}
 N_{11}^{-}  & \cdots & N_{1m}^{-}  \\
 \vdots & \cdots & \vdots \\
 N_{m1}^{-}  & \cdots & N_{mm}^{-} 
 \end{array}
 \right)
 \end{equation}
 where each $N_{ik}^{\pm}$ is an $\# C_i \times \# C_k$-matrix with nonnegative integer entries whose  row  sum is $n^{\pm}_{ik}$.   \\
 (ii) The activation and repression matrices $W^{+}$ and $W^{-}$ have block structures 
 $$
 W^{+} = 
 \left( 
 \begin{array}{ccc}
 Q_{11}^{+}  & \cdots & Q_{1m}^{+}  \\
 \vdots & \cdots & \vdots \\
 Q_{m1}^{+}  & \cdots & Q_{mm}^{+} 
 \end{array}
 \right) ,\qquad 
 W^{-} = 
 \left( 
 \begin{array}{ccc}
 Q_{11}^{-}  & \cdots & Q_{1m}^{-}  \\
 \vdots & \cdots & \vdots \\
 Q_{m1}^{-}  & \cdots & Q_{mm}^{-} 
 \end{array}
 \right)
 $$
 where the matrices $Q^{+}_{ij}$,  $Q^{-}_{ij}$ have nonnegative real entries satisfying the following: for 
  $i=1, \ldots, m$,  the product of the {nonzero} entries of each row of $\left( Q^{+}_{i1} \ldots  Q^{+}_{im} {Q^{-}_{i1} \ldots Q^{-}_{im}}\right)$ equals  the product  of the {nonzero}  entries of 
 $\left( q^{+}_{i1} \ldots q^{+}_{im} \, q^{-}_{i1} \ldots q^{-}_{im}\right)$.\\
 \end{thm}

\begin{proof}
Direct application of Proposition~\ref{prop:MULTsyn} and Theorem 2.13 in \cite{AD18}. 
\end{proof}

\begin{exam} 
{Let $Q$ be the 3}-gene GRN  for the MULT model equation (considered in Example~\ref{ex:quo_not_prod}) with the $3 \times 3$ multiplicity matrices:
$$N^+ = 
\left(
\begin{array}{r|r|r}
2 & 0 &  0\\
\hline
0 & 0 & 0 \\
\hline 
0 & 0 & 0 
\end{array}
\right), \  
N^{-} = 
\left(
\begin{array}{r|r|r}
0 & 1 & 1\\
\hline
1 & 0 & 0 \\
\hline 
1 & 0 & 0
\end{array}
\right) 
$$
and  activation and repression weighted adjacency matrices  
$$
Q^+ = 
\left(
\begin{array}{r|r|r}
1& 0 &  0\\
\hline
0  & 0 & 0 \\
\hline 
0  & 0 & 0 
\end{array}
\right), \  
Q^{-} = 
\left(
\begin{array}{r|r|r}
 0 & 3 & 3\\
\hline
2 & 0 & 0 \\
\hline 
3 & 0 & 0
\end{array}
\right)\, .$$
{Considering the 4-gene GRNs that are lifts of $Q$ for the PROD (MULT) model, we have that 
 any 4-}gene GRN with $4 \times 4$ multiplicity matrices 
$$
M^+ = 
\left(
\begin{array}{rr|r|r}
1& 1& 0 &  0\\
1 & 1 & 0 & 0 \\
\hline
0 & 0 & 0 & 0 \\
\hline 
0 & 0 & 0 & 0 
\end{array}
\right), \quad 
M^{-} = 
\left(
\begin{array}{rr|r|r}
0& 0 & 1 & 1\\
0 & 0 & 1 & 1 \\
\hline
0 & 1 & 0 & 0 \\
\hline 
1 & 0 & 0 & 0
\end{array}
\right)
$$
and activation and repression weighted adjacency matrices 
$$
W^+ = 
\left(
\begin{array}{rr|r|r}
w_{11}^{+} & w_{12}^{+} & 0 &  0\\
w_{21}^{+} & w_{22}^{+} & 0 & 0 \\
\hline
0 & 0 & 0 & 0 \\
\hline 
0 & 0 & 0 & 0 
\end{array}
\right), \quad 
W^{-} = 
\left(
\begin{array}{rr|r|r}
0& 0 & w_{13}^{-} & w_{14}^{-}\\
0 & 0 & w_{23}^{-} & w_{24}^{-}\\
\hline
0 & 2 & 0 & 0 \\
\hline 
3 & 0 & 0 & 0
\end{array}
\right), 
$$
where $w_{11}^{+} w_{12}^{+} w_{13}^{-} w_{14}^{-} = w_{21}^{+} w_{22}^{+} w_{23}^{-} w_{24}^{-} =9$, is a lift of the 3-gene GRN. By Proposition~\ref{prop:MULTsyn}, assuming that genes $1,2$ have the same internal dynamics,  any such lift has the synchronization partition $P$ with the parts $\{1,2\},\, \{3\},\, \{ 4\},\, \{5\}$ and the restriction of MULT model equations  to the synchrony space $\Delta_P = \{ x:\, x_1 = x_2\}$ gives rise to the MULT model equations (\ref{eq:notPROD}). 
\hfill $\Diamond$
\end{exam}

\section{Regulatory-dependent synchrony spaces}\label{sec:regdepsyn} 

Gene regulation functions provide an unique characteristic of GRNs relating concentrations of transcription factors such as activators or repressors to the promoter activities. They have been investigated in two ways (\cite{T_2010}). Classical molecular biology explores mechanistic details of transcription and translation activities (\cite{H_2007,P_2004}), whereas the emerging field of system biology quantifies gene expressions on a larger scale in a statistical mechanical framework without requiring physical details of macromolecular interactions (\cite{BT_2004,SRG_2008,Numbers_2005, Numbers_2005b}).

So far, assuming that the regulator functions $\mathrm{act}$ and $\mathrm{rep}$  are Hill-like according to  (\ref{eq;reg_Hill_like_0})-(\ref{eq;reg_Hill_like_infty}), we have shown,
for the GRN equations models SUM and MULT (includes PROD),  that 
the synchrony spaces are completely determined by their structural matrices (adjacency matrices $W^\pm$ in case of  SUM and PROD models; multiplicity matrices $M^\pm$ and adjacency matrices $W^\pm$ in case of MULT model). See Propositions \ref{prop:synSUM} and \ref{prop:MULTsyn}. 

For regulator functions that are not Hill-like this may not be true:
there may be other synchrony subspaces forced by those regulator functions, which we call {\it regulatory-dependent synchrony spaces}. 

We illustrate this phenomena with the next example.

\begin{figure}[!htb]
		\centerline{
			\includegraphics[width=.8\textwidth]{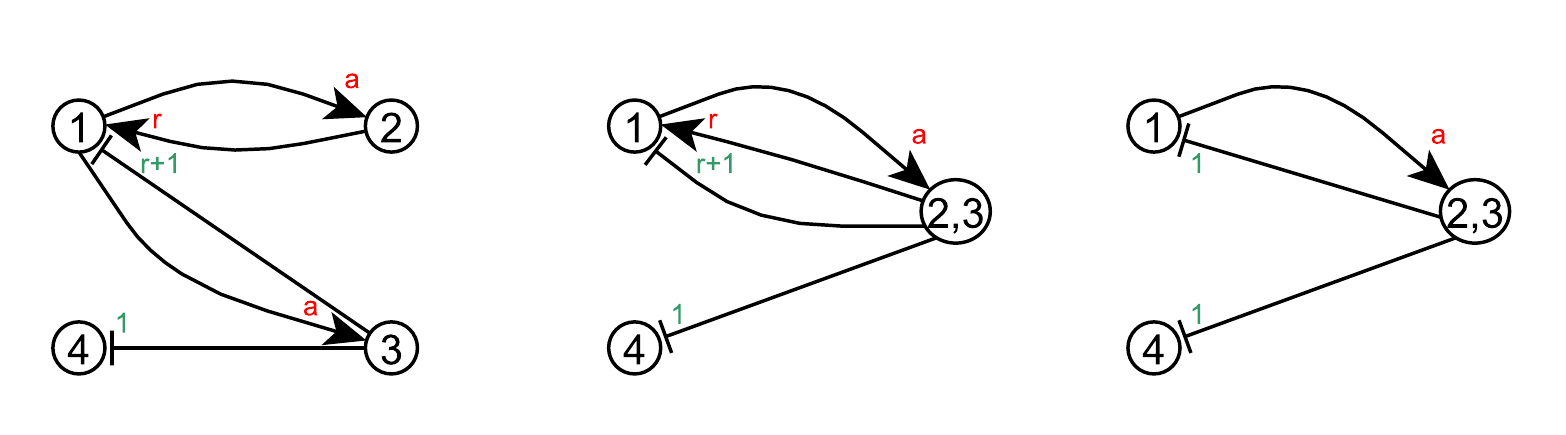}
			}
		\caption{Left: The graph representation of the 4-gene MULT model GRN equations (\ref{eq:ex_exotic}). 
		When restricted to $S_o=\{x_2=x_3\}$, the system  (\ref{eq:ex_exotic})  has two different graph representations: one corresponds to  (\ref{eq:ex_exotic_c})  
		(middle) and the other corresponds to (\ref{eq:ex_exotic_b}) (right). This is due to the particular relation between the activation and repression functions that are considered in equations (\ref{eq:ex_exotic}).}
		\label{F:exotic}  
	\end{figure}

\begin{exam}\label{ex:exotic} 
Consider the repressor and activation functions given by (\ref{eq:5_rep})-(\ref{eq:5_act}), in Section~\ref{sec:mammalian} of the GRN model for the mammalian circadian oscillator,
and the GRN MULT model equations
	\begin{equation} \label{eq:ex_exotic}
	\begin{cases}
	\dot x_1=\act(x_2,r)\rep(x_3,r) -x_1\\
	\dot x_2=\act(x_1,a)-x_2\\
	\dot x_3=\act(x_1,a)-x_3\\
	\dot x_4=\rep(x_3,1)-x_4
	\end{cases}
	\end{equation}
	for any $a\ge 1$, $r\ge 0$ 
	which can be represented by a multi-arrowed graph. See Figure \ref{F:exotic}(left). 	
	We have that
$$S_o=\{x_2=x_3\},\qquad S=\{x_1=x_4,x_2=x_3\}$$ 
	are synchrony subspaces of (\ref{eq:ex_exotic}). More precisely, $S_o$ is a synchrony subspace, since $P = \left\{ \{1\}, \{2, 3\}, \{4\}\right\}$ is a synchronization partition for the MULT model equations (\ref{eq:ex_exotic}). On the other hand, $S$ is flow-invariant for the equations (\ref{eq:ex_exotic}) due to the fact that ${\mathrm act} (x,r) {\mathrm rep} (x,r) = {\mathrm rep}(x,1)$ holds for all $r,x\ge 0$.   Thus, $S$ is a regulatory dependent synchrony pattern for equations (\ref{eq:ex_exotic}). Note that $I_1^+=\{2\}$ while $I_4^+=\emptyset$. 	The synchrony subspace $S_o=\{x_2=x_3\}$ can be recognized on the graph, since cells $2,3$ have identical positive input sets.
	However,  $S=\{x_1=x_4,x_2=x_3\}$ may not be immediately identified {looking at the} graph.
	
	One way to observe that $S$ is a synchrony subspace, is to start with $S_o=\{x_2=x_3\}$ and consider its quotient network, which can be represented by two different graphs. Indeed, when restricted to $S_o$, (\ref{eq:ex_exotic}) becomes
	\begin{equation} \label{eq:ex_exotic_c}
	\begin{cases}
	\dot x_1=\act(x_{2,3},r)\rep(x_{2,3},r) -x_1\\
	\dot x_{2,3}=\act(x_1,a)-x_{2,3}\\
	\dot x_4=\rep(x_{2,3},1)-x_4
	\end{cases}
	\end{equation}
	or equivalently,
	\begin{equation} \label{eq:ex_exotic_b}
	\begin{cases}
	\dot x_1=\rep(x_{2,3},1) -x_1\\
	\dot x_{2,3}=\act(x_1,a)-x_{2,3}\\
	\dot x_4=\rep(x_{2,3},1)-x_4
	\end{cases}
	\end{equation}
	due to the relation $\act(x,r)\rep(x,r)=\rep(x,1)$. They correspond to Figure \ref{F:exotic} middle and right, respectively. 
	\hfill $\Diamond$	
\end{exam} 

As shown by Example \ref{ex:exotic}, regulatory dependent synchrony spaces(which are not inferred by the network structure) can exist, which are directly related to the modelling of the network itself including the choice of regulator functions and choice of models. Also,  
as shown in Example \ref{ex:funny},  hidden relations between $\act$ and $\rep$ functions can cause  cancellations of regulator functions, which results in the ``de-coupling'' from regulating genes.

\section{Conclusions and  Outlook} \label{sec:future} 

	This work is a first contribution to show how the theory of coupled cell systems can play a relevant role in the study of GRNs. 
We considered two dynamical models of GRNs, depending on whether the gene regulatory is additive (the SUM model) or multiplicative (the MULT model).   Using theoretical results from coupled cell networks, we analyse the robust patterns of synchrony supported by these gene regulatory models and found out that the gene synchronization patterns can be quite different for the SUM and MULT models in general. Moreover, we have shown that other unexpected synchrony patterns may occur  when the activation and regression functions satisfy specific relations.

Related to the process of gene duplication and the phenomena of genetic redundancy, we explored the concept of quotient networks and network lifting in 
 both  SUM and MULT models.  From our results, it follows in particular that, if a SUM or MULT model of equations for a small GRN presents, for example, oscillatory behaviour, then  we can enumerate bigger GRNs that admit this small GRN as a quotient network while preserving the initial oscillatory behaviour with some of the genes being synchronized. That is, the bigger GRN is a lifting of the smaller one where it is guaranteed that the bigger GRN has a synchrony space such that the associated quotient SUM or MULT equations restricted to the synchrony space are precisely the SUM or MULT model equations describing the small GRN presenting the oscillatory behaviour. Furthermore, we described a method for constructing such network lifts. It was pointed out that there is no uniqueness in this process and so the choice of the lifts can be adapted to the particular types of applications under consideration. 

Although our results are theoretical and by no means indicate how they stand in relation to empirical evidence or experimental data in GRNs,  we wish to remark that coupled cell network formalisms and GRNs do share common features in structural dynamical properties. Therefore, it could be worthwhile incorporating theoretical considerations using coupled cell networks in  the discussion of dynamical processes that have been addressed in GRNs.

In a future work, motivated by the work presented in \cite{KKBLRH_2014} (see Remark~\ref{rem:mot_osc}), we aim to explore other potential approaches of embedding an oscillatory gene network with fewer genes into a bigger gene network, without using synchrony subspaces and where all genes may oscillate in a non-synchronized fashion.

\vspace{5mm}

\noindent {\bf Acknowledgments} \\
The authors MA and AD were partially supported by CMUP (UID/MAT/00144/2013), which is funded by FCT (Portugal) with national (MEC) and European structural funds (FEDER), under the partnership agreement PT2020. HR would like to thank  Dr. Sabine Le Borne at Technical University of Hamburg, Germany for supporting her independent research and wish to express her gratitude to Hefei University, China for dynamic research exchanges and communications.

\end{document}